\documentclass[12pt]{article}
\usepackage{amsmath}
\usepackage{amsfonts}
\usepackage{amssymb}

\def\thebibliograph#1#2{\section*{{\normalsize \bf #2}}\list
   {[\arabic{enumi}]}{\settowidth\labelwidth{[#1]}\leftmargin\labelwidth
     \advance\leftmargin\labelsep
     \usecounter{enumi}}
     \def\newblock{\hskip .11em plus .33em minus -.07em}
     \sloppy
     \sfcode`\.=1000\relax}

\newtheorem{theorem}{Theorem}

\newtheorem{definition}{Definition}
\newtheorem{corollary}{Corollary}
\newtheorem{lemma}{Lemma}
\newtheorem{remark}{Remark}
\input{tcilatex}
\begin{document}

\title{Besov-type spaces with variable smoothness and integrability II}
\author{ Douadi Drihem \ \thanks{%
M'sila University, Department of Mathematics, Laboratory of Functional
Analysis and Geometry of Spaces , P.O. Box 166, M'sila 28000, Algeria,
e-mail: \texttt{\ douadidr@yahoo.fr}}}
\date{\today }
\maketitle

\begin{abstract}
The aim of this paper is to study properties of Besov-type spaces with
variable smoothness. We show that these spaces are characterized by the $%
\varphi $-transforms in appropriate sequence spaces and we obtain atomic
decompositions for these spaces.
\end{abstract}

\section{Introduction}

The most known general scales of function spaces are the scales of Besov
spaces and Triebel-Lizorkin spaces and it is known that they cover many
well-known classical function spaces such as H\"{o}lder-Zygmund spaces and
Sobolev spaces, see Triebel's monographes \cite{T1} and \cite{T2} for the
history of these function spaces. These spaces play an important role in
Harmonic Analysis.

The theory of these spaces had a remarkable development in part due to its
usefulness in applications. For instance, they appear in the study of
partial differential equations.

\ \ In recent years, there has been growing interest in generalizing
classical spaces such as Lebesgue, Sobolev spaces, Besov spaces and
Triebel-Lizorkin spaces to the case with either variable integrability or
variable smoothness. The motivation for the increasing interest in such
spaces comes not only from theoretical purposes, but also from applications
to fluid dynamics \cite{Ru00}, image restoration and PDE with non-standard
growth conditions.

Variable Besov-type spaces $B_{p(\cdot ),q(\cdot )}^{\alpha (\cdot ),\tau
(\cdot )}$ have been introduced in \cite{D5}, where their basic properties
are given, such as the Sobolev type embeddings and \ that under some
conditions these spaces are just the Besov spaces $B_{\infty ,\infty
}^{\alpha (\cdot )+n(1/\tau (\cdot )-1/p(\cdot ))}$. For constant exponents,
these spaces unify and generalize many classical function spaces including
Besov spaces, Besov-Morrey spaces (see, for example, \cite[Corollary 3.3]%
{WYY}).

\ \ The main aim of this paper is to present another essential property of
the Besov-type spaces with variable smoothness and integrability such as the 
$\varphi $-transforms characterization and the atomic decomposition.\vskip5pt

\ \ The paper is organized as follows. First we give some preliminaries
where we fix some notations and recall some basics facts on function spaces
with variable integrability\ and we give some key technical lemmas needed in
the proofs of the main statements. For making the presentation clearer, we
give their proofs later in Section 5. We then define the Besov-type spaces $%
B_{p(\cdot ),q(\cdot )}^{\alpha (\cdot ),\tau (\cdot )}$ and $\widetilde{B}%
_{p(\cdot ),q(\cdot )}^{\alpha (\cdot ),p(\cdot )}$ as follows. Let $%
\mathcal{Q}$ be the set of all dyadic cubes in $\mathbb{R}^{n}$. For each
cube $P\in \mathcal{Q}$, we denote by $\chi _{P}$ the characteristic
function of $P$. Select a pair of Schwartz functions $\Phi $ and $\varphi $
satisfy%
\begin{equation}
\text{supp}\mathcal{F}\Phi \subset \overline{B(0,2)}\text{ and }|\mathcal{F}%
\Phi (\xi )|\geq c\text{ if }|\xi |\leq \frac{5}{3}  \label{Ass1}
\end{equation}%
and 
\begin{equation}
\text{supp}\mathcal{F}\varphi \subset \overline{B(0,2)}\backslash B(0,1/2)%
\text{ and }|\mathcal{F}\varphi (\xi )|\geq c\text{ if }\frac{3}{5}\leq |\xi
|\leq \frac{5}{3}  \label{Ass2}
\end{equation}%
where $c>0$. Let $\alpha :\mathbb{R}^{n}\rightarrow \mathbb{R}$, $p,q,\tau
\in \mathcal{P}_{0}$ and $\Phi $ and $\varphi $ satisfy $\mathrm{\eqref{Ass1}%
}$ and $\mathrm{\eqref{Ass2}}$, respectively and we put $\varphi
_{v}=2^{vn}\varphi (2^{v}\cdot )$.\textrm{\ }The Besov-type space $%
\widetilde{B}_{p(\cdot ),q(\cdot )}^{\alpha (\cdot ),p(\cdot )}$\ is the
collection of all $f\in \mathcal{S}^{\prime }(\mathbb{R}^{n})$\ such that 
\begin{equation*}
\left\Vert f\right\Vert _{\widetilde{B}_{p(\cdot ),q(\cdot )}^{\alpha (\cdot
),p(\cdot )}}:=\sup_{P\in \mathcal{Q}}\left\Vert \left( \frac{2^{v\alpha
\left( \cdot \right) }\varphi _{v}\ast f}{|P|^{1/p(\cdot )}}\chi _{P}\right)
_{v\geq v_{P}^{+}}\right\Vert _{\ell ^{q(\cdot )}(L^{p(\cdot )})}<\infty ,
\end{equation*}%
where $\varphi _{0}$ is replaced by $\Phi $. While the Besov-type space $%
B_{p(\cdot ),q(\cdot )}^{\alpha (\cdot ),\tau (\cdot )}$\ is the collection
of all $f\in \mathcal{S}^{\prime }(\mathbb{R}^{n})$\ such that 
\begin{equation*}
\left\Vert f\right\Vert _{B_{p(\cdot ),q(\cdot )}^{\alpha (\cdot ),\tau
(\cdot )}}:=\sup_{P\in \mathcal{Q}}\left\Vert \left( \frac{2^{v\alpha \left(
\cdot \right) }\varphi _{v}\ast f}{\left\Vert \chi _{P}\right\Vert _{\tau
(\cdot )}}\chi _{P}\right) _{v\geq v_{P}^{+}}\right\Vert _{\ell ^{q(\cdot
)}(L^{p(\cdot )})}<\infty ,
\end{equation*}%
where $\varphi _{0}$ is replaced by $\Phi $, see Section 2 for the
definition of $\mathcal{P}_{0}$, $\ell ^{q(\cdot )}(L^{p(\cdot )})$ and $%
v_{P}^{+}$. In this section several basic properties such as the $\varphi $%
-transform characterization are obtained. The main statements are formulated
in Section 4, where we give the atomic decomposition of these function
spaces. It is shown that the element $f\in \mathcal{S}^{\prime }(\mathbb{R}%
^{n})$ in the space $B_{p(\cdot ),q(\cdot )}^{\alpha (\cdot ),\tau (\cdot )}$
or $\widetilde{B}_{p(\cdot ),q(\cdot )}^{\alpha (\cdot ),p(\cdot )}$can be
represented as%
\begin{equation*}
f=\sum\limits_{v=0}^{\infty }\sum\limits_{m\in \mathbb{Z}^{n}}\lambda
_{v,m}\varrho _{v,m},\text{ \ \ \ \ converging in }\mathcal{S}^{\prime }(%
\mathbb{R}^{n})\text{,}
\end{equation*}%
where $\varrho _{v,m}$'s are the so-called atoms and the sequence complex
numbers $\{\lambda _{v,m}\}$ belongs to an appropriate sequence space.
Moreover, based on these sequence spaces equivalent quasi-norms for
corresponding function spaces are derived. In this section we also give some
key technical lemmas needed in the proofs of the main statements.

\section{Preliminaries}

As usual, we denote by $\mathbb{R}^{n}$ the $n$-dimensional real Euclidean
space, $\mathbb{N}$ the collection of all natural numbers and $\mathbb{N}%
_{0}=\mathbb{N}\cup \{0\}$. The letter $\mathbb{Z}$ stands for the set of
all integer numbers.\ The expression $f\lesssim g$ means that $f\leq c\,g$
for some independent constant $c$ (and non-negative functions $f$ and $g$),
and $f\approx g$ means $f\lesssim g\lesssim f$. As usual for any $x\in 
\mathbb{R}$, $[x]$ stands for the largest integer smaller than or equal to $%
x $.\vskip5pt

For $x\in \mathbb{R}^{n}$ and $r>0$ we denote by $B(x,r)$ the open ball in $%
\mathbb{R}^{n}$ with center $x$ and radius $r$. By supp $f$ we denote the
support of the function $f$ , i.e., the closure of its non-zero set. If $%
E\subset {\mathbb{R}^{n}}$ is a measurable set, then $|E|$ stands for the
(Lebesgue) measure of $E$ and $\chi _{E}$ denotes its characteristic
function.\vskip5pt

The symbol $\mathcal{S}(\mathbb{R}^{n})$ is used in place of the set of all
Schwartz functions $\phi $ on $\mathbb{R}^{n}$, i.e., $\phi $ is infinitely
differentiable and%
\begin{equation*}
\left\Vert \phi \right\Vert _{\mathcal{S}_{M}}=\sup_{\gamma \in \mathbb{N}%
_{0}^{n},|\gamma |\leq M}\sup_{x\in \mathbb{R}^{n}}|\partial ^{\gamma }\phi
(x)|(1+|x|)^{n+M+|\gamma |}<\infty
\end{equation*}%
for all $M\in \mathbb{N}$. We denote by $\mathcal{S}^{\prime }(\mathbb{R}%
^{n})$ the dual space of all tempered distributions on $\mathbb{R}^{n}$. We
define the Fourier transform of a function $f\in \mathcal{S}(\mathbb{R}^{n})$
by $\mathcal{F}(f)(\xi )=(2\pi )^{-n/2}\int_{\mathbb{R}^{n}}e^{-ix\cdot \xi
}f(x)dx$. Its inverse is denoted by $\mathcal{F}^{-1}f$. Both $\mathcal{F}$
and $\mathcal{F}^{-1}$ are extended to the dual Schwartz space $\mathcal{S}%
^{\prime }(\mathbb{R}^{n})$ in the usual way.\vskip5pt

The Hardy-Littlewood maximal operator $\mathcal{M}$ is defined on $L_{%
\mathrm{loc}}^{1}$ by%
\begin{equation*}
\mathcal{M}f(x)=\sup_{r>0}\frac{1}{\left\vert B(x,r)\right\vert }%
\int_{B(x,r)}\left\vert f(y)\right\vert dy.
\end{equation*}%
For $v\in \mathbb{Z}$ and $m=(m_{1},...,m_{n})\in \mathbb{Z}^{n}$, let $%
Q_{v,m}$ be the dyadic cube in $\mathbb{R}^{n}$, $Q_{v,m}=%
\{(x_{1},...,x_{n}):m_{i}\leq 2^{v}x_{i}<m_{i}+1,i=1,2,...,n\}$. For the
collection of all such cubes we use $\mathcal{Q}:=\{Q_{v,m}:v\in \mathbb{Z}%
,m\in \mathbb{Z}^{n}\}$. For each cube $Q$, we denote by $x_{Q_{v,m}}$ the
lower left-corner $2^{-v}m$ of $Q=Q_{v,m}$, its side length by $l(Q)$ and
for $r>0$, we denote by $rQ$ the cube concentric with $Q$ having the side
length $rl(Q)$. Furthermore, we put $v_{Q}=-\log _{2}l(Q)$ and $%
v_{Q}^{+}=\max (v_{Q},0)$.\vskip5pt

For $v\in \mathbb{Z}$, $\varphi \in \mathcal{S}(\mathbb{R}^{n})$ and $x\in 
\mathbb{R}^{n}$, we set $\widetilde{\varphi }(x):=\overline{\varphi (-x)}$, $%
\varphi _{v}(x):=2^{vn}\varphi (2^{v}x)$, and%
\begin{equation*}
\varphi _{v,m}(x):=2^{vn/2}\varphi (2^{v}x-m)=|Q_{v,m}|^{1/2}\varphi
_{v}(x-x_{Q_{v,m}})\quad \text{if\quad }Q=Q_{v,m}.
\end{equation*}

By $c$ we denote generic positive constants, which may have different values
at different occurrences. Although the exact values of the constants are
usually irrelevant for our purposes, sometimes we emphasize their dependence
on certain parameters (e.g. $c(p)$ means that $c$ depends on $p$, etc.).
Further notation will be properly introduced whenever needed.

The variable exponents that we consider are always measurable functions $p$
on $\mathbb{R}^{n}$ with range in $[c,\infty \lbrack $ for some $c>0$. We
denote the set of such functions by $\mathcal{P}_{0}$. The subset of
variable exponents with range $[1,\infty \lbrack $ is denoted by $\mathcal{P}
$. We use the standard notation $p^{-}:=\underset{x\in \mathbb{R}^{n}}{\text{%
ess-inf}}$ $p(x)$,$\quad p^{+}:=\underset{x\in \mathbb{R}^{n}}{\text{ess-sup 
}}p(x)$.

The variable exponent modular is defined by $\varrho _{p(\cdot )}(f):=\int_{%
\mathbb{R}^{n}}\rho _{p(x)}(\left\vert f(x)\right\vert )dx$, where $\rho
_{p}(t)=t^{p}$. The variable exponent Lebesgue space $L^{p(\cdot )}$\
consists of measurable functions $f$ on $\mathbb{R}^{n}$ such that $\varrho
_{p(\cdot )}(\lambda f)<\infty $ for some $\lambda >0$. We define the
Luxemburg (quasi)-norm on this space by the formula $\left\Vert f\right\Vert
_{p(\cdot )}:=\inf \Big\{\lambda >0:\varrho _{p(\cdot )}\Big(\frac{f}{%
\lambda }\Big)\leq 1\Big\}$. A useful property is that $\left\Vert
f\right\Vert _{p(\cdot )}\leq 1$ if and only if $\varrho _{p(\cdot )}(f)\leq
1$, see \cite{DHHR}, Lemma 3.2.4.

Let $p,q\in \mathcal{P}_{0}$. The mixed Lebesgue-sequence space $\ell
^{q(\cdot )}(L^{p(\cdot )})$ is defined on sequences of $L^{p(\cdot )}$%
-functions by the modular%
\begin{equation*}
\varrho _{\ell ^{q(\cdot )}(L^{p\left( \cdot \right)
})}((f_{v})_{v}):=\sum\limits_{v}\inf \Big\{\lambda _{v}>0:\varrho _{p(\cdot
)}\Big(\frac{f_{v}}{\lambda _{v}^{1/q(\cdot )}}\Big)\leq 1\Big\}.
\end{equation*}%
The (quasi)-norm is defined from this as usual:%
\begin{equation}
\left\Vert \left( f_{v}\right) _{v}\right\Vert _{\ell ^{q(\cdot
)}(L^{p\left( \cdot \right) })}:=\inf \Big\{\mu >0:\varrho _{\ell ^{q(\cdot
)}(L^{p(\cdot )})}\Big(\frac{1}{\mu }(f_{v})_{v}\Big)\leq 1\Big\}.
\label{mixed-norm}
\end{equation}%
If $q^{+}<\infty $, then we can replace $\mathrm{\eqref{mixed-norm}}$ by the
simpler expression $\varrho _{\ell ^{q(\cdot )}(L^{p(\cdot
)})}((f_{v})_{v}):=\sum\limits_{v}\left\Vert |f_{v}|^{q(\cdot )}\right\Vert
_{\frac{p(\cdot )}{q(\cdot )}}$. Furthermore, if $p$ and $q$ are constants,
then $\ell ^{q(\cdot )}(L^{p(\cdot )})=\ell ^{q}(L^{p})$. The case $%
p:=\infty $ can be included by replacing the last modular by $\varrho _{\ell
^{q(\cdot )}(L^{\infty })}((f_{v})_{v}):=\sum\limits_{v}\left\Vert
\left\vert f_{v}\right\vert ^{q(\cdot )}\right\Vert _{\infty }$.

It is known, cf. \cite{AH} and \cite{KV121}, that $\ell ^{q(\cdot
)}(L^{p(\cdot )})$ is a norm if $q(\cdot )\geq 1$ is constant almost
everywhere (a.e.) on $\mathbb{R}^{n}$ and $p(\cdot )\geq 1$, or if $\frac{1}{%
p(x)}+\frac{1}{q(x)}\leq 1$ a.e. on $\mathbb{R}^{n}$, or if $1\leq q(x)\leq
p(x)<\infty $ a.e. on $\mathbb{R}^{n}$.

We say that $g:\mathbb{R}^{n}\rightarrow \mathbb{R}$ is \textit{locally }log%
\textit{-H\"{o}lder continuous}, abbreviated $g\in C_{\text{loc}}^{\log }$,
if there exists $c_{\log }(g)>0$ such that%
\begin{equation}
\left\vert g(x)-g(y)\right\vert \leq \frac{c_{\log }(g)}{\log
(e+1/\left\vert x-y\right\vert )}  \label{lo-log-Holder}
\end{equation}%
for all $x,y\in \mathbb{R}^{n}$. We say that $g$ satisfies the log\textit{-H%
\"{o}lder decay condition}, if there exists $g_{\infty }\in \mathbb{R}$ and
a constant $c_{\log }>0$ such that%
\begin{equation*}
\left\vert g(x)-g_{\infty }\right\vert \leq \frac{c_{\log }}{\log
(e+\left\vert x\right\vert )}
\end{equation*}%
for all $x\in \mathbb{R}^{n}$. We say that $g$ is \textit{globally}-log%
\textit{-H\"{o}lder continuous}, abbreviated $g\in C^{\log }$, if it is%
\textit{\ }locally log-H\"{o}lder continuous and satisfies the log-H\"{o}%
lder decay\textit{\ }condition.\textit{\ }The constants $c_{\log }(g)$ and $%
c_{\log }$ are called the \textit{locally }log\textit{-H\"{o}lder constant }%
and the log\textit{-H\"{o}lder decay constant}, respectively\textit{.} We
note that all functions $g\in C_{\text{loc}}^{\log }$ always belong to $%
L^{\infty }$.\vskip5pt

We define the following class of variable exponents%
\begin{equation*}
\mathcal{P}^{\mathrm{log}}:=\Big\{p\in \mathcal{P}:\frac{1}{p}\text{ is
globally-log-H\"{o}lder continuous}\Big\},
\end{equation*}%
were introduced in $\mathrm{\cite[Section \ 2]{DHHMS}}$. We define $%
1/p_{\infty }:=\lim_{|x|\rightarrow \infty }1/p(x)$\ and we use the
convention $\frac{1}{\infty }=0$. Note that although $\frac{1}{p}$ is
bounded, the variable exponent $p$ itself can be unbounded. It was shown in $%
\mathrm{\cite{DHHR}}$\textrm{, }Theorem 4.3.8 that $\mathcal{M}:L^{p(\cdot
)}\rightarrow L^{p(\cdot )}$ is bounded if $p\in \mathcal{P}^{\mathrm{log}}$
and $p^{-}>1$, see also $\mathrm{\cite{DHHMS}}$, Theorem 1.2.\ Also if $p\in 
\mathcal{P}^{\mathrm{log}}$, then the convolution with a radially decreasing 
$L^{1}$-function is bounded on $L^{p(\cdot )}$: $\Vert \varphi \ast f\Vert _{%
{p(\cdot )}}\leq c\Vert \varphi \Vert _{{1}}\Vert f\Vert _{{p(\cdot )}}$. We
also refer to the papers $\mathrm{\cite{CFMP}}$ and $\mathrm{\cite{Di}}$%
\textrm{,} where various results on maximal function in variable Lebesgue
spaces were obtained.\vskip5pt

It is known that for $p\in \mathcal{P}^{\mathrm{log}}$ we have%
\begin{equation}
\Vert \chi _{B}\Vert _{{p(\cdot )}}\Vert \chi _{B}\Vert _{{p}^{\prime }{%
(\cdot )}}\approx |B|.  \label{DHHR}
\end{equation}%
Also,%
\begin{equation}
\Vert \chi _{B}\Vert _{{p(\cdot )}}\approx |B|^{\frac{1}{p(x)}},\quad x\in B
\label{DHHR1}
\end{equation}%
for small balls $B\subset {\mathbb{R}^{n}}$ ($|B|\leq 2^{n}$), and 
\begin{equation}
\Vert \chi _{B}\Vert _{{p(\cdot )}}\approx |B|^{\frac{1}{p_{\infty }}}
\label{DHHR2}
\end{equation}%
for large balls ($|B|\geq 1$), with constants only depending on the $\log $-H%
\"{o}lder constant of $p$ (see, for example, \cite[Section 4.5]{DHHR}). Here 
${p}^{\prime }$ denotes the conjugate exponent of $p$ given by $1/{p(\cdot )}%
+1/{p}^{\prime }{(\cdot )}=1$. These properties are hold if $p\in \mathcal{P}%
_{0}^{\mathrm{log}}$, since $\Vert \chi _{B}\Vert _{{p(\cdot )}}=\Vert \chi
_{B}\Vert _{{p(\cdot )/a}}^{1/a}$ and $\frac{p}{a}{\in }\mathcal{P}^{\mathrm{%
log}}$ if $p^{-}\geq a.$\vskip5pt

Recall that $\eta _{v,m}(x):=2^{nv}(1+2^{v}\left\vert x\right\vert )^{-m}$,
for any $x\in \mathbb{R}^{n}$, $v\in \mathbb{N}_{0}$ and $m>0$. Note that $%
\eta _{v,m}\in L^{1}$ when $m>n$ and that $\left\Vert \eta _{v,m}\right\Vert
_{1}=c_{m}$ is independent of $v$, where this type of function was
introduced in \cite{HN07} and \cite{DHHR}.

\subsection{Some technical lemmas}

In this subsection we present some results which are useful for us. The
following lemma is from \cite[Lemma 19]{KV122}, see also \cite[Lemma 6.1]%
{DHR}.

\begin{lemma}
\label{DHR-lemma}Let $\alpha \in C_{\mathrm{loc}}^{\log }$ and let $R\geq
c_{\log }(\alpha )$, where $c_{\log }(\alpha )$ is the constant from $%
\mathrm{\eqref{lo-log-Holder}}$ for $\alpha $. Then%
\begin{equation*}
2^{v\alpha (x)}\eta _{v,m+R}(x-y)\leq c\text{ }2^{v\alpha (y)}\eta
_{v,m}(x-y)
\end{equation*}%
with $c>0$ independent of $x,y\in \mathbb{R}^{n}$ and $v,m\in \mathbb{N}%
_{0}. $
\end{lemma}

The previous lemma allows us to treat the variable smoothness in many cases
as if it were not variable at all, namely we can move the term inside the
convolution as follows:%
\begin{equation*}
2^{v\alpha (x)}\eta _{v,m+R}\ast f(x)\leq c\text{ }\eta _{v,m}\ast
(2^{v\alpha (\cdot )}f)(x).
\end{equation*}

\begin{lemma}
\label{r-trick}Let $r,R,N>0$, $m>n$ and $\theta ,\omega \in \mathcal{S}%
\left( \mathbb{R}^{n}\right) $ with $\mathrm{supp}\mathcal{F}\omega \subset 
\overline{B(0,1)}$. Then there exists $c=c(r,m,n)>0$ such that for all $g\in 
\mathcal{S}^{\prime }\left( \mathbb{R}^{n}\right) $, we have%
\begin{equation}
\left\vert \theta _{R}\ast \omega _{N}\ast g\left( x\right) \right\vert \leq
c\text{ }\max \Big(1,\Big(\frac{N}{R}\Big)^{m}\Big)(\eta _{N,m}\ast
\left\vert \omega _{N}\ast g\right\vert ^{r}(x))^{1/r},\quad x\in \mathbb{R}%
^{n},  \label{r-trick-est}
\end{equation}%
where $\theta _{R}(\cdot )=R^{n}\theta (R\cdot )$, $\omega _{N}(\cdot
)=N^{n}\omega (N\cdot )$ and $\eta _{N,m}:=N^{n}(1+N\left\vert \cdot
\right\vert )^{-m}$.
\end{lemma}

This lemma is a slight variant of \cite[Chapter V, Theorem 5]{ST89}, see
also \cite[Lemma A.7]{DHR}. For the convenience of the reader, we give the
proof in the Appendix.

The following lemma is from \cite[Lemma 2.11]{D5}.

\begin{lemma}
\label{n-cube-est}Let $\mathbb{\tau }\in \mathcal{P}_{0}^{\log }$ and $k\in 
\mathbb{Z}^{n}$.

$\mathrm{(i)}$ For any cubes $P$ and $Q$, we have%
\begin{equation*}
\frac{\left\Vert \chi _{P+kl(Q)}\right\Vert _{\tau (\cdot )}}{\left\Vert
\chi _{P}\right\Vert _{\tau (\cdot )}}\leq c\left( 1+\frac{l(Q)}{l(P)}%
|k|\right) ^{c_{\log }(\frac{1}{\tau} )}
\end{equation*}%
with $c>0$ independent of $l(Q)$, $l(P)$ and $k$.

$\mathrm{(ii)}$ For any cubes $P$ and $Q$, such that $P\subset Q$, we have%
\begin{equation*}
C\left( \frac{|Q|}{|P|}\right) ^{1/\tau ^{+}}\leq \frac{\left\Vert \chi
_{Q}\right\Vert _{\tau (\cdot )}}{\left\Vert \chi _{P}\right\Vert _{\tau
(\cdot )}}\leq c\left( \frac{|Q|}{|P|}\right) ^{1/\tau ^{-}}
\end{equation*}%
with $c,C>0$ are independent of $|Q|$ and $|P|$.
\end{lemma}

Let $L_{\tau (\cdot )}^{p(\cdot )}$ be the collection of functions $f\in L_{%
\text{loc}}^{p(\cdot )}(\mathbb{R}^{n})$ such that%
\begin{equation*}
\left\Vert f\right\Vert _{L_{\tau (\cdot )}^{p(\cdot )}}:=\sup \left\Vert 
\frac{f\chi _{P}}{\left\Vert \chi _{P}\right\Vert _{\tau (\cdot )}}%
\right\Vert _{p(\cdot )}<\infty ,\quad p,\tau \in \mathcal{P}_{0},
\end{equation*}%
where the supremum is taken over all dyadic cubes $P$ with $|P|\geq 1$.
Also, the spaces $\widetilde{L^{p(\cdot )}}$ is defined to be the set of all
function $f$ such that%
\begin{equation*}
\left\Vert f\right\Vert _{\widetilde{L^{p(\cdot )}}}:=\sup \left\Vert f\chi
_{P}\right\Vert _{p(\cdot )}<\infty ,\quad p\in \mathcal{P}_{0},
\end{equation*}%
where the supremum is taken over all dyadic cubes $P$ with $|P|=1$. Notice
that 
\begin{equation}
\left\Vert f\right\Vert _{L_{\tau (\cdot )}^{p(\cdot )}}\leq
1\Leftrightarrow \sup_{P\in \mathcal{Q},|P|\geq 1}\left\Vert \left\vert 
\frac{f}{\left\Vert \chi _{P}\right\Vert _{\tau (\cdot )}}\right\vert
^{q(\cdot )}\chi _{P}\right\Vert _{p(\cdot )/q(\cdot )}\leq 1.
\label{mod-est}
\end{equation}

Let $\theta _{R}$ be as in Lemma \ref{r-trick}.

\begin{lemma}
\label{key-estimate1}Let $R,N>0$, $\mathbb{\tau },p\in \mathcal{P}_{0}^{\log
}$, $0<r<p^{-}$ and $\theta ,\omega \in \mathcal{S}(\mathbb{R}^{n})$ with $%
\mathrm{supp}\mathcal{F}\omega \subset \overline{B(0,1)}$.

$\mathrm{(i)}$ For any $f\in \mathcal{S}^{\prime }(\mathbb{R}^{n})$, any $%
m>2n+c_{\log }(\frac{1}{\tau })r$ and any dyadic cube $P$ with $|P|\geq 1$,
we have%
\begin{equation*}
\left\Vert \frac{\theta _{R}\ast \omega _{N}\ast f}{\left\Vert \chi
_{P}\right\Vert _{\tau (\cdot )}}\chi _{P}\right\Vert _{p(\cdot )}\leq c%
\text{ }\max \Big(1,\Big(\frac{N}{R}\Big)^{m}\Big)\max (1,\left(
Nl(P)\right) ^{(n-m)/r})\left\Vert \omega _{N}\ast f\right\Vert _{L_{\tau
(\cdot )}^{p(\cdot )}},
\end{equation*}%
such that the right-hand side is finite, where $c>0$ is independent of $R,$ $%
N$ and $l(P).$

$\mathrm{(ii)}$ For any $f\in \mathcal{S}^{\prime }(\mathbb{R}^{n})$, any $%
m>2n$ and any dyadic cube $P$ with $|P|=1$, we have%
\begin{equation*}
\left\Vert (\theta _{R}\ast \omega _{N}\ast f)\chi _{P}\right\Vert _{p(\cdot
)}\leq c\text{ }\max \Big(1,\Big(\frac{N}{R}\Big)^{m}\Big)\max
(1,N^{(n-m)/r})\left\Vert \omega _{N}\ast f\right\Vert _{\widetilde{%
L^{p(\cdot )}}},
\end{equation*}%
such that the right-hand side is finite, where $c>0$ is independent of $R$
and $N$.
\end{lemma}

The proof of this lemma is postponed to the Appendix. We introduce the
abbreviations%
\begin{equation*}
\left\Vert \left( f_{v}\right) _{v}\right\Vert _{\ell ^{q(\cdot
)}(L_{p(\cdot )}^{p(\cdot )})}:=\sup_{\{P\in \mathcal{Q},|P|\leq
1\}}\left\Vert \left( \frac{f_{v}}{|P|^{1/p(\cdot )}}\chi _{P}\right)
_{v\geq v_{P}}\right\Vert _{\ell ^{q(\cdot )}(L^{p(\cdot )})}
\end{equation*}%
and%
\begin{equation*}
\left\Vert \left( f_{v}\right) _{v}\right\Vert _{\ell ^{\tau (\cdot
),q(\cdot )}(L^{p(\cdot )})}:=\sup_{P\in \mathcal{Q}}\left\Vert \left( \frac{%
f_{v}}{\left\Vert \chi _{P}\right\Vert _{\tau (\cdot )}}\chi _{P}\right)
_{v\geq v_{P}^{+}}\right\Vert _{\ell ^{q(\cdot )}(L^{p(\cdot )})},
\end{equation*}%
where, $v_{P}=-\log _{2}l(P)$ and $v_{P}^{+}=\max (v_{P},0)$. The following
lemma is the $\ell ^{q(\cdot )}(L_{p(\cdot )}^{p(\cdot )})$(-$\ell ^{\tau
(\cdot ),q(\cdot )}(L^{p(\cdot )})$)-version of Lemma 4.7 from A. Almeida
and P. H\"{a}st\"{o} \cite{AH} (we use it, since the maximal operator is in
general not bounded on $\ell ^{q(\cdot )}(L^{p(\cdot )})$, see \cite[Example
4.1]{AH}).

\begin{lemma}
\label{Alm-Hastolemma1}Let $p\in \mathcal{P}^{\log }$, $q,\tau \in \mathcal{P%
}_{0}^{\log }$ with $0<q^{-}\leq q^{+}<\infty $ and $p^{-}>1$.\newline
$\mathrm{(i)}$ For $m>2n+c_{\log }(1/\tau )+c_{\log }(1/q)$, there exists $%
c>0$ such that%
\begin{equation*}
\left\Vert (\eta _{v,m}\ast f_{v})_{v}\right\Vert _{\ell ^{\tau (\cdot
),q(\cdot )}(L^{p(\cdot )})}\leq c\left\Vert (f_{v})_{v}\right\Vert _{\ell
^{\tau (\cdot ),q(\cdot )}(L^{p(\cdot )})}.
\end{equation*}%
\newline
$\mathrm{(ii)}$ For $m>2n+c_{\log }(1/p)+c_{\log }(1/q)$, there exists $c>0$
such that%
\begin{equation*}
\left\Vert (\eta _{v,m}\ast f_{v})_{v}\right\Vert _{\ell ^{q(\cdot
)}(L_{p(\cdot )}^{p(\cdot )})}\leq c\left\Vert (f_{v})_{v}\right\Vert _{\ell
^{q(\cdot )}(L_{p(\cdot )}^{p(\cdot )})}.
\end{equation*}
\end{lemma}

The proof (i) is given in \cite[Lemma 2.12]{D5}, their arguments are true to
prove (ii) in view of the fact that $\left\Vert \chi _{P}\right\Vert
_{p(\cdot )}\approx |P|^{1/p(\cdot )}$, since the supremum taken with
respect to dyadic cubes with side length $\leq 1$.

The next three lemmas are from \cite{DHR} where the first tells us that in
most circumstances two convolutions are as good as one.

\begin{lemma}
\label{Conv-est}For $v_{0},v_{1}\in \mathbb{N}_{0}$ and $m>n$, we have%
\begin{equation*}
\eta _{v_{0},m}\ast \eta _{v_{1},m}\approx \eta _{\min (v_{0},v_{1}),m}
\end{equation*}%
with the constant depending only on $m$ and $n$.
\end{lemma}

\begin{lemma}
\label{Conv-est1}Let $v\in\mathbb{N}_{0}$ and $m>n$. Then for any $Q\in 
\mathcal{Q}$ with $l(Q) =2^{-v}$, $y\in Q$ and $x\in \mathbb{R}^{n}$, we have%
\begin{equation*}
\eta _{v,m}\ast \left( \frac{\chi _{Q}}{|Q|}\right) (x)\approx \eta
_{v,m}(x-y)
\end{equation*}%
with the constant depending only on $m$ and $n$.
\end{lemma}

\begin{lemma}
\label{Conv-est2}Let $v,j\in \mathbb{N}_{0}$, $r\in (0,1]$ and $m>\frac{n}{r}
$. Then for any $Q\in \mathcal{Q}$ with $l(Q)=2^{-v}$, we have%
\begin{equation*}
(\eta _{j,m}\ast \eta _{v,m}\ast \chi _{Q})^{r}\approx
2^{(v-j)^{+}n(1-r)}\eta _{j,mr}\ast \eta _{v,mr}\ast \chi _{Q},
\end{equation*}%
where the constant depends only on $m$, $n$ and $r$.
\end{lemma}

The next lemma is a Hardy-type inequality which is easy to prove.

\begin{lemma}
\label{lq-inequality}\textit{Let }$0<a<1,J\in \mathbb{Z}$\textit{\ and }$%
0<q\leq \infty $\textit{. Let }$\left\{ \varepsilon _{k}\right\} $\textit{\
be a sequences of positive real numbers and denote} $\delta
_{k}=\sum_{j=J^{+}}^{k}a^{k-j}\varepsilon _{j}$, $k\geq J^{+}$.\textit{\ }%
Then there exists constant $c>0\ $\textit{depending only on }$a$\textit{\
and }$q$ such that%
\begin{equation*}
\left( \sum\limits_{k=J^{+}}^{\infty }\delta _{k}^{q}\right) ^{1/q}\leq c%
\text{ }\left( \sum\limits_{k=J^{+}}^{\infty }\varepsilon _{k}^{q}\right)
^{1/q}.
\end{equation*}
\end{lemma}

\begin{lemma}
\label{Key-lemma}Let $\alpha \in C_{\mathrm{loc}}^{\log }$ and $p,q,\mathbb{%
\tau }\in \mathcal{P}_{0}^{\log }$ with $0<q^{-}\leq q^{+}<\infty $. Let $%
\left\{ f_{k}\right\} _{k\in \mathbb{N}_{0}}$ be a sequence of measurable
functions on $\mathbb{R}^{n}$. For all $v\in \mathbb{N}_{0}$ and $x\in 
\mathbb{R}^{n}$, let $g_{v}(x)=\sum_{k=0}^{\infty }2^{-|k-v|\delta }f_{k}(x)$%
. Then there exists a positive constant $c$, independent of $\left\{
f_{k}\right\} _{k\in \mathbb{N}_{0}}$ such that%
\begin{equation*}
\left\Vert (g_{v})_{v}\right\Vert _{\ell ^{\tau (\cdot ),q(\cdot
)}(L^{p(\cdot )})}\leq c\left\Vert (f_{v})_{v}\right\Vert _{\ell ^{\tau
(\cdot ),q(\cdot )}(L^{p(\cdot )})},\quad \delta >\frac{n}{\tau ^{-}}.
\end{equation*}%
and 
\begin{equation*}
\left\Vert (g_{v})_{v}\right\Vert _{\ell ^{q(\cdot )}(L_{p(\cdot )}^{p(\cdot
)})}\leq c\left\Vert (f_{v})_{v}\right\Vert _{\ell ^{q(\cdot )}(L_{p(\cdot
)}^{p(\cdot )})},\quad \delta >\frac{n}{p^{-}}.
\end{equation*}
\end{lemma}

The proof of Lemma \ref{Key-lemma} is postponed to the Appendix.

\section{The spaces\textbf{\ }$\widetilde{B}_{p(\cdot ),q(\cdot )}^{\protect%
\alpha (\cdot ),p(\cdot )}$ and $B_{p(\cdot ),q(\cdot )}^{\protect\alpha %
(\cdot ),\protect\tau (\cdot )}$}

In this section we\ present the Fourier analytical definition of Besov-type
spaces of variable smoothness and integrability\ and we prove the basic
properties in analogy to the Besov-type spaces with fixed exponents. Select
a pair of Schwartz functions $\Phi $ and $\varphi $ satisfy $\mathrm{%
\eqref{Ass1}}$ and $\mathrm{\eqref{Ass2}}$, respectively. It easy to see
that $\int_{\mathbb{R}^{n}}x^{\gamma }\varphi (x)dx=0$ for all multi-indices 
$\gamma \in \mathbb{N}_{0}^{n}$. For the convenience of the reader we repeat
the definition of the spaces\textbf{\ }$\widetilde{B}_{p(\cdot ),q(\cdot
)}^{\alpha (\cdot ),p(\cdot )}$ and $B_{p(\cdot ),q(\cdot )}^{\alpha (\cdot
),\tau (\cdot )}$.

\begin{definition}
\label{B-F-def}Let $\alpha :\mathbb{R}^{n}\rightarrow \mathbb{R}$, $p,q,\tau
\in \mathcal{P}_{0}$ and $\Phi $ and $\varphi $ satisfy $\mathrm{\eqref{Ass1}%
}$ and $\mathrm{\eqref{Ass2}}$, respectively and we put $\varphi
_{v}=2^{vn}\varphi (2^{v}\cdot )$.

$\mathrm{(i)}$ The Besov-type space $\widetilde{B}_{p(\cdot ),q(\cdot
)}^{\alpha (\cdot ),p(\cdot )}$\ is the collection of all $f\in \mathcal{S}%
^{\prime }(\mathbb{R}^{n})$\ such that 
\begin{equation}
\left\Vert f\right\Vert _{\widetilde{B}_{p(\cdot ),q(\cdot )}^{\alpha (\cdot
),p(\cdot )}}:=\sup_{P\in \mathcal{Q}}\left\Vert \left( \frac{2^{v\alpha
\left( \cdot \right) }\varphi _{v}\ast f}{|P|^{1/p(\cdot )}}\chi _{P}\right)
_{v\geq v_{P}^{+}}\right\Vert _{\ell ^{q(\cdot )}(L^{p(\cdot )})}<\infty ,
\label{B-def}
\end{equation}%
where $\varphi _{0}$ is replaced by $\Phi $.

$\mathrm{(ii)}$ The Besov-type space $B_{p(\cdot ),q(\cdot )}^{\alpha (\cdot
),\tau (\cdot )}$\ is the collection of all $f\in \mathcal{S}^{\prime }(%
\mathbb{R}^{n})$\ such that 
\begin{equation*}
\left\Vert f\right\Vert _{B_{p(\cdot ),q(\cdot )}^{\alpha (\cdot ),\tau
(\cdot )}}:=\sup_{P\in \mathcal{Q}}\left\Vert \left( \frac{2^{v\alpha \left(
\cdot \right) }\varphi _{v}\ast f}{\left\Vert \chi _{P}\right\Vert _{\tau
(\cdot )}}\chi _{P}\right) _{v\geq v_{P}^{+}}\right\Vert _{\ell ^{q(\cdot
)}(L^{p(\cdot )})}<\infty ,
\end{equation*}%
where $\varphi _{0}$ is replaced by $\Phi $.
\end{definition}

Using the system $\{\varphi _{v}\}_{v\in \mathbb{N}_{0}}$ we can define the
norm%
\begin{equation*}
\left\Vert f\right\Vert _{B_{p,q}^{\alpha ,\tau }}:=\sup_{P\in \mathcal{Q}}%
\frac{1}{\left\vert P\right\vert ^{\tau }}\left(
\sum\limits_{v=v_{P}^{+}}^{\infty }2^{v\alpha q}\left\Vert \left( \varphi
_{v}\ast f\right) \chi _{P}\right\Vert _{p}^{q}\right) ^{1/q}
\end{equation*}%
for constants $\alpha $ and $p,q\in (0,\infty ]$. The Besov-type space $%
B_{p,q}^{\alpha ,\tau }$ consist of all distributions $f\in \mathcal{S}%
^{\prime }(\mathbb{R}^{n})$ for which $\left\Vert f\right\Vert
_{B_{p,q}^{\alpha ,\tau }}<\infty $. It is well-known that these spaces do
not depend on the choice of the system $\{\varphi _{v}\}_{v\in \mathbb{N}%
_{0}}$ (up to equivalence of quasinorms). Further details on the classical
theory of these spaces can be found in \cite{D1} and \cite{WYY}; see also 
\cite{D4} for recent developments.

One recognizes immediately that if $\alpha $, $\tau $, $p$ and $q$ are
constants, then $\widetilde{B}_{p(\cdot ),q(\cdot )}^{\alpha (\cdot
),p(\cdot )}=B_{p,q}^{\alpha ,1/p}$ and $B_{p(\cdot ),q(\cdot )}^{\alpha
(\cdot ),\tau (\cdot )}=B_{p,q}^{\alpha ,\tau }$. When, $q:=\infty $\ the
Besov-type space $\widetilde{B}_{p(\cdot ),\infty }^{\alpha (\cdot ),p(\cdot
)}$\ consist of all distributions $f\in \mathcal{S}^{\prime }(\mathbb{R}%
^{n}) $\ such that 
\begin{equation*}
\sup_{P\in \mathcal{Q},v\geq v_{P}^{+}}\left\Vert \frac{2^{v\alpha \left(
\cdot \right) }\varphi _{v}\ast f}{|P|^{1/p(\cdot )}}\chi _{P}\right\Vert
_{p(\cdot )}<\infty
\end{equation*}%
and the Besov-type space $B_{p(\cdot ),\infty }^{\alpha (\cdot ),\tau (\cdot
)}$\ consist of all distributions $f\in \mathcal{S}^{\prime }(\mathbb{R}%
^{n}) $\ such that 
\begin{equation*}
\sup_{P\in \mathcal{Q},v\geq v_{P}^{+}}\left\Vert \frac{2^{v\alpha \left(
\cdot \right) }\varphi _{v}\ast f}{\left\Vert \chi _{P}\right\Vert _{\tau
(\cdot )}}\chi _{P}\right\Vert _{p(\cdot )}<\infty .
\end{equation*}%
Let $B_{J}$ be any ball of $\mathbb{R}^{n}$ with radius $2^{-J}$, $J\in 
\mathbb{Z}$. In the definition of the spaces $B_{p(\cdot ),q(\cdot
)}^{\alpha (\cdot ),\tau (\cdot )}$ and $\widetilde{B}_{p(\cdot ),q(\cdot
)}^{\alpha (\cdot ),p(\cdot )}$ if we replace the dyadic cubes $P$ by the
balls $B_{J}$, then we obtain equivalent quasi-norms. From these if we
replace dyadic cubes $P$ in Definition \ref{B-F-def} by arbitrary cubes $P$,
we then obtain equivalent quasi-norms.

The spaces $B_{p(\cdot ),q(\cdot )}^{\alpha (\cdot ),\tau (\cdot )}$, were
introduced and studied in \cite{D5}, where we proved that our spaces are
well-defined, i.e., independent of the choice of the resolution of unity\
and we gave some properties of these function spaces, see Theorem \ref%
{new-est1}, below. Moreover the Sobolev embeddings for these function spaces
are obtained. While the first time we introduce the spaces $\widetilde{B}%
_{p(\cdot ),q(\cdot )}^{\alpha (\cdot ),p(\cdot )}$ with the quasi-norm $%
\mathrm{\eqref{B-def}}$. Independently, D. Yang, C. Zhuo and W. Yuan, \cite%
{YZW15} studied the function spaces $B_{p(\cdot ),q(\cdot )}^{\alpha (\cdot
),\tau (\cdot )}$\thinspace\ where several properties are obtained such as
atomic decomposition and the boundedness of trace operator.

Moreover, the following remarkable features are given in \cite{D5} where
these results with fixed exponents are given in \cite{YY13} and \cite{WYY}.

\begin{theorem}
\label{new-est1}Let $\alpha \in C_{\mathrm{loc}}^{\log }$, $p,p_{1},p_{2},q,%
\mathbb{\tau }\in \mathcal{P}_{0}^{\log }$ and $0<q^{+}<\infty $.\newline
$\mathrm{(i)}$ Let $\mathbb{\tau }_{\infty }\in (0,p^{-}]$. If $(1/\tau
-1/p)^{-}>0$ or $(1/\tau -1/p)^{-}\geq 0$ and $q:=\infty $, then $B_{p(\cdot
),q(\cdot )}^{\alpha (\cdot ),\tau (\cdot )}=B_{\infty ,\infty }^{\alpha
(\cdot )+n(1/\tau (\cdot )-1/p(\cdot ))}$, with equivalent norms.\newline
$\mathrm{(ii)}$ If $(p_{2}-p_{1})^{+}\leq 0$, then $B_{p_{2}(\cdot ),q(\cdot
)}^{\alpha (\cdot )+n/\tau (\cdot )+n/p_{2}(\cdot )-n/p_{1}(\cdot
)}\hookrightarrow B_{p_{1}(\cdot ),q(\cdot )}^{\alpha (\cdot ),\tau (\cdot
)}.$ \newline
$\mathrm{(iii)}$ We have%
\begin{equation*}
B_{p(\cdot ),q(\cdot )}^{\alpha (\cdot ),\tau (\cdot )}\hookrightarrow
B_{\infty ,\infty }^{\alpha (\cdot )+n/\tau (\cdot )-n/p(\cdot )}.
\end{equation*}
\end{theorem}

Here $B_{p(\cdot ),q(\cdot )}^{\alpha (\cdot )}$ is the Besov space of
variable smoothness and integrability and it is the collection of all $f\in 
\mathcal{S}^{\prime }(\mathbb{R}^{n})$\ such that 
\begin{equation*}
\left\Vert f\right\Vert _{B_{p(\cdot ),q(\cdot )}^{\alpha (\cdot
)}}:=\left\Vert \left( 2^{v\alpha \left( \cdot \right) }\varphi _{v}\ast
f\right) _{v\geq 0}\right\Vert _{\ell ^{q(\cdot )}(L^{p(\cdot )})}<\infty ,
\end{equation*}%
which introduced and investigated in \cite{AH} and \cite{KV122}\ for further
results. Taking $\alpha \in \mathbb{R}$ and $q\in (0,\infty ]$ as constants
we derive the spaces $B_{p(\cdot ),q}^{\alpha }$ studied by Xu in \cite{Xu08}
and \cite{Xu09}. We refer the reader to the recent paper \cite{YHSY} for
further details, historical remarks and more references on embeddings of
Besov-type spaces with fixed exponents.

Let $0<u\leq p<\infty $. The Morrey space $\mathcal{M}_{u}^{p}$ is defined
to be the set of all $u$-locally Lebesgue-integrable functions $f$ on $%
\mathbb{R}^{n}$ such that%
\begin{equation*}
\left\Vert f\right\Vert _{\mathcal{M}_{u}^{p}}:=\sup_{B}|B|^{\frac{1}{p}-%
\frac{1}{u}}\left( \int_{B}|f(x)|^{u}dx\right) ^{1/u}<\infty ,
\end{equation*}%
where the supremum is taken over all balls $B$ in $\mathbb{R}^{n}$. The
spaces $\mathcal{M}_{u}^{p}$ are quasi-Banach spaces (Banach spaces for $%
u\geq 1$). They were introduced by Morrey in \cite{M38} and belong to the
wider class of Morrey-Campanato spaces, cf. \cite{P}. They can be considered
as a complement to $L^{p}$ spaces. As a matter of fact, $\mathcal{M}%
_{p}^{p}=L^{p}$. One can easily see that $\mathcal{M}_{w}^{p}\hookrightarrow 
\mathcal{M}_{u}^{p}\ $if\ $0<u\leq w<\infty .$

\begin{definition}
Let $\alpha :\mathbb{R}^{n}\rightarrow \mathbb{R},0<u\leq p<\infty $ and $%
0<q\leq \infty $. Let $\Phi $ and $\varphi $ satisfy $\mathrm{\eqref{Ass1}}$
and $\mathrm{\eqref{Ass2}}$, respectively and we put $\varphi
_{v}=2^{vn}\varphi (2^{v}\cdot )$. The Besov-Morrey space $\mathcal{N}%
_{p,q,u}^{\alpha (\cdot )}$\ is the collection of all $f\in \mathcal{S}%
^{\prime }(\mathbb{R}^{n})$\ such that 
\begin{equation*}
\left\Vert f\right\Vert _{\mathcal{N}_{p,q,u}^{\alpha (\cdot )}}:=\left(
\sum_{v=0}^{\infty }\left\Vert 2^{v\alpha \left( \cdot \right) }\varphi
_{v}\ast f\right\Vert _{\mathcal{M}_{u}^{p}}^{q}\right) ^{1/q}<\infty ,
\end{equation*}
where $\varphi _{0}$ is replaced by $\Phi $.
\end{definition}

Besov-Morrey spaces with fixed exponents were introduced by Netrusov \cite%
{N84}. Kozono and Yamazaki \cite{KY95} studied semilinear heat equations and
Navier-Stokes equations with initial data belonging to Besov-Morrey spaces.
The investigations were continued by Mazzucato \cite{Ma03}, where one can
find the wavelet decomposition of Besov-Morrey spaces. On the other hand,
the Besov-Morrey space $\mathcal{N}_{p,q,u}^{\alpha }$ is a proper subspace
of the space $B_{u,q}^{\alpha ,\frac{1}{u}-\frac{1}{p}}$ with $u<p$, and $%
q<\infty $, see \cite{SYY10}. Further properties for these function spaces
can be found in \cite{Sa08}, \cite{Sa09} and \cite{Sa10}.

Recently, Triebel in \cite{T4} further introduced and studied some local
versions of these smoothness Morrey-type spaces and also considered their
applications in heat equations and Navier-Stokes equations. More recent
results can be found in \cite{YSY13}, where they studied the relations
between Triebel's local spaces and the Besov-type and Triebel-Lizorkin-type
spaces and their associated uniform spaces.

D. Yang and W. Yuan introduced and investigated in \cite{YY1} and \cite{YY2}
the homogeneous Besov and Triebel-Lizorkin spaces, which generalize the
homogeneous Besov and Triebel-Lizorkin spaces.

The Besov-Morrey spaces with variable exponents have been first introduced
in \cite{FX11}, where are introduced equivalent quasi-norms of these new
spaces, which are formulated in terms of Peetre's maximal functions. Also
the authors obtain the atomic, molecular and wavelet decompositions of these
new spaces.

In the next proposition we present the relations between variable
Besov-Morrey spaces and variable Besov-type spaces, see \cite{D5}.

\begin{proposition}
Let $\alpha \in C_{\mathrm{loc}}^{\log }$, $0<q<\infty $ and $0<p<u<\infty .$%
\newline
$\mathrm{(i)}$ For $0<q<\infty $ we have the continuous embeddings%
\begin{equation*}
\mathcal{N}_{u,q,p}^{\alpha (\cdot )}\hookrightarrow B_{p,q}^{\alpha
(\cdot ),\left( \frac{1}{p}-\frac{1}{u}\right) ^{-1}}\text{.}
\end{equation*}%
$\mathrm{(ii)}$ We have
\begin{equation*}
\mathcal{N}_{u,\infty ,p}^{\alpha (\cdot )}=B_{p,\infty }^{\alpha
(\cdot ),\left( \frac{1}{p}-\frac{1}{u}\right) ^{-1}}.
\end{equation*}
\end{proposition}

Sometimes it is of great service if one can restrict sup$_{P\in \mathcal{Q}}$
in the definition to a supremum taken with respect to dyadic cubes with side
length $\leq 1$.

\begin{lemma}
\label{new-equinorm copy(1)}Let $\alpha \in C_{\mathrm{loc}}^{\log }$ and $%
p,q,\mathbb{\tau }\in \mathcal{P}_{0}^{\log }$ with $\mathbb{\tau }_{\infty
}\in (0,p^{-}]$ and $0<q^{+}<\infty .$

$\mathrm{(i)}$ A tempered distribution $f$ belongs to $\widetilde{B}%
_{p(\cdot ),q(\cdot )}^{\alpha (\cdot ),p(\cdot )}$ if and only if,%
\begin{equation*}
\left\Vert f\right\Vert _{\widetilde{B}_{p(\cdot ),q(\cdot )}^{\alpha (\cdot
),p(\cdot )}}^{\#}:=\sup_{\{P\in \mathcal{Q},|P|\leq 1\}}\left\Vert \left( 
\frac{2^{v\alpha \left( \cdot \right) }\varphi _{v}\ast f}{|P|^{1/p(\cdot )}}%
\chi _{P}\right) _{v\geq v_{P}}\right\Vert _{\ell ^{q(\cdot )}(L^{p(\cdot
)})}<\infty .
\end{equation*}%
Furthermore, the quasi-norms $\left\Vert f\right\Vert _{\widetilde{B}%
_{p(\cdot ),q(\cdot )}^{\alpha (\cdot ),p(\cdot )}}$ and $\left\Vert
f\right\Vert _{\widetilde{B}_{p(\cdot ),q(\cdot )}^{\alpha (\cdot ),p(\cdot
)}}^{\#}$ are equivalent.

$\mathrm{(ii)}$ A tempered distribution $f$ belongs to $B_{p(\cdot ),q(\cdot
)}^{\alpha (\cdot ),\tau (\cdot )}$ if and only if,%
\begin{equation*}
\left\Vert f\right\Vert _{B_{p(\cdot ),q(\cdot )}^{\alpha (\cdot ),\tau
(\cdot )}}^{\#}:=\sup_{\{P\in \mathcal{Q},|P|\leq 1\}}\left\Vert \left( 
\frac{2^{v\alpha \left( \cdot \right) }\varphi _{v}\ast f}{\left\Vert \chi
_{P}\right\Vert _{\tau (\cdot )}}\chi _{P}\right) _{v\geq v_{P}}\right\Vert
_{\ell ^{q(\cdot )}(L^{p(\cdot )})}<\infty .
\end{equation*}%
Furthermore, the quasi-norms $\left\Vert f\right\Vert _{B_{p(\cdot ),q(\cdot
)}^{\alpha (\cdot ),\tau (\cdot )}}$ and $\left\Vert f\right\Vert
_{B_{p(\cdot ),q(\cdot )}^{\alpha (\cdot ),\tau (\cdot )}}^{\#}$ are
equivalent.
\end{lemma}

The proof is similar to that of \cite{D5}. We omit the details.

\begin{remark}
$\mathrm{(i)}$ We like to point out that this result with fixed exponents is
given in \cite[Lemma 2.2]{WYY} with $1/\tau $ in place of $\tau $.

$\mathrm{(ii)}$ Let $\alpha \in C_{\mathrm{loc}}^{\log }$, $p,q,\in \mathcal{%
P}_{0}^{\log }$ and $0<q^{+}<\infty $. As in \cite{D5}, we obtain$\ 
\widetilde{B}_{p(\cdot ),\infty }^{\alpha (\cdot ),p(\cdot )}=B_{\infty
,\infty }^{\alpha (\cdot )}$. Also, $2^{v(\alpha (x)+n(1/\tau
(x)-1/p(x)))}|\varphi _{v}\ast f(x)|\leq c\left\Vert f\right\Vert
_{B_{p(\cdot ),q(\cdot )}^{\alpha (\cdot ),\tau (\cdot )}}$ for any $x\in 
\mathbb{R}^{n}$, $\alpha \in C_{\mathrm{loc}}^{\log }$ and $p,q\in \mathcal{P%
}_{0}^{\log }.$

$\mathrm{(iii)}$ It is clear that if $\alpha $ and $p$ are constants, then $%
\widetilde{B}_{p(\cdot ),p(\cdot )}^{\alpha (\cdot ),p(\cdot )}=F_{\infty
,p}^{\alpha }$, see \cite{FJ90} for the properties of $F_{\infty ,p}^{\alpha
}$.

$\mathrm{(iv)}$ We can easily prove that if $\alpha \in C_{\mathrm{loc}%
}^{\log }$, $p,q\in \mathcal{P}_{0}^{\log }$ and $0<q^{+}<\infty $, then $%
B_{p(\cdot ),q(\cdot )}^{\alpha (\cdot ),p(\cdot )}\hookrightarrow 
\widetilde{B}_{p(\cdot ),q(\cdot )}^{\alpha (\cdot ),p(\cdot )}$.

$\mathrm{(v)}$ In \cite{D5} the definition of Besov-type spaces $B_{p(\cdot
),q(\cdot )}^{\alpha (\cdot ),\tau (\cdot )}$ is based on the technique of
decomposition of unity.
\end{remark}

Let $\alpha \in C_{\mathrm{loc}}^{\log },p,q\in \mathcal{P}_{0}^{\log }$ and 
$\alpha _{0}<{\alpha }^{-}$. We obtain%
\begin{equation*}
\widetilde{B}_{{p(\cdot )},q{(\cdot )}}^{{\alpha (\cdot ),p(\cdot )}%
}\hookrightarrow \widetilde{B}_{{p(\cdot )},\infty }^{{\alpha }_{0}{,p(\cdot
)}}=B_{{\infty },\infty }^{{\alpha }_{0}}\hookrightarrow \mathcal{S}^{\prime
}(\mathbb{R}^{n}).
\end{equation*}%
Let ${\alpha }^{+}<\alpha _{1}$. We obtain%
\begin{equation*}
\mathcal{S}(\mathbb{R}^{n})\hookrightarrow B_{{\infty },\infty }^{{\alpha }%
_{1}}=\widetilde{B}_{{p(\cdot )},\infty }^{{\alpha }_{1}{,p(\cdot )}%
}\hookrightarrow \widetilde{B}_{{p(\cdot )},q{(\cdot )}}^{{\alpha (\cdot
),p(\cdot )}}.
\end{equation*}%
We use $A_{p(\cdot ),q(\cdot )}^{\alpha (\cdot ),\tau (\cdot )}$ to denote
either $\widetilde{B}_{p(\cdot ),q(\cdot )}^{\alpha (\cdot ),p(\cdot )}$ or $%
B_{p(\cdot ),q(\cdot )}^{\alpha (\cdot ),\tau (\cdot )}$.

\begin{theorem}
Let $\alpha \in C_{\mathrm{loc}}^{\log }$ and $p,q,\mathbb{\tau }\in 
\mathcal{P}_{0}^{\log }$ with $0<q^{+}<\infty $. Then 
\begin{equation*}
\mathcal{S}(\mathbb{R}^{n})\hookrightarrow A_{{p(\cdot )},q{(\cdot )}}^{{%
\alpha (\cdot ),\tau (\cdot )}}\hookrightarrow \mathcal{S}^{\prime }(\mathbb{%
R}^{n}).
\end{equation*}
\end{theorem}

Similar arguments of \cite{D5} can be used to prove the following
Sobolev-type embeddings.

\begin{theorem}
Let $\alpha _{0},\alpha _{1}\in C_{\mathrm{loc}}^{\log }$ and $%
p_{0},p_{1},q\in \mathcal{P}_{0}^{\log }$ with $0<q^{+}<\infty $. If ${%
\alpha }_{0}>{\alpha }_{1}$\ and ${\alpha }_{0}{(x)-}\frac{n}{p_{0}(x)}={%
\alpha }_{1}{(x)-}\frac{n}{p_{1}(x)}$ with $\Big(\frac{p_{0}}{p_{1}}\Big)%
^{-}<1$, then 
\begin{equation*}
\widetilde{B}_{{p}_{0}{(\cdot )},q{(\cdot )}}^{{\alpha }_{0}{(\cdot
),p(\cdot )}}\hookrightarrow \widetilde{B}_{{p}_{1}{(\cdot )},q{(\cdot )}}^{{%
\alpha }_{1}{(\cdot ),p(\cdot )}}.
\end{equation*}
\end{theorem}

Notice that the case of $B_{p(\cdot ),q(\cdot )}^{\alpha (\cdot ),\tau
(\cdot )}$ spaces is given in \cite{D5}.

Let $\Phi $ and $\varphi $ satisfy, respectively $\mathrm{\eqref{Ass1}}$ and 
$\mathrm{\eqref{Ass2}}$. By \cite[pp. 130--131]{FJ90}, there exist \
functions $\Psi \in \mathcal{S}(\mathbb{R}^{n})$ satisfying $\mathrm{%
\eqref{Ass1}}$ and $\psi \in \mathcal{S}(\mathbb{R}^{n})$ satisfying $%
\mathrm{\eqref{Ass2}}$ such that for all $\xi \in \mathbb{R}^{n}$%
\begin{equation}
\mathcal{F}\widetilde{\Phi }(\xi )\mathcal{F}\Psi (\xi )+\sum_{j=1}^{\infty }%
\mathcal{F}\widetilde{\varphi }(2^{-j}\xi )\mathcal{F}\psi (2^{-j}\xi
)=1,\quad \xi \in \mathbb{R}^{n}.  \label{Ass4}
\end{equation}

Furthermore, we have the following identity for all $f\in \mathcal{S}%
^{\prime }(\mathbb{R}^{n})$; see \cite[(12.4)]{FJ90}%
\begin{eqnarray*}
f &=&\Psi \ast \widetilde{\Phi }\ast f+\sum_{v=1}^{\infty }\psi _{v}\ast 
\widetilde{\varphi }_{v}\ast f \\
&=&\sum_{m\in \mathbb{Z}^{n}}\widetilde{\Phi }\ast f(m)\Psi (\cdot
-m)+\sum_{v=1}^{\infty }2^{-vn}\sum_{m\in \mathbb{Z}^{n}}\widetilde{\varphi }%
_{v}\ast f(2^{-v}m)\psi _{v}(\cdot -2^{-v}m).
\end{eqnarray*}%
Recall that the $\varphi $-transform $S_{\varphi }$ is defined by setting $%
(S_{\varphi })_{0,m}=\langle f,\Phi _{m}\rangle $ where $\Phi _{m}(x)=\Phi
(x-m)$ and $(S_{\varphi })_{v,m}=\langle f,\varphi _{v,m}\rangle $ where $%
\varphi _{v,m}(x)=2^{vn/2}\varphi (2^{v}x-m)$ and $v\in \mathbb{N}$. The
inverse $\varphi $-transform $T_{\psi }$ is defined by 
\begin{equation*}
T_{\psi }\lambda =\sum_{m\in \mathbb{Z}^{n}}\lambda _{0,m}\Psi
_{m}+\sum_{v=1}^{\infty }\sum_{m\in \mathbb{Z}^{n}}\lambda _{v,m}\psi _{v,m},
\end{equation*}%
where $\lambda =\{\lambda _{v,m}\in \mathbb{C}:v\in \mathbb{N}_{0},m\in 
\mathbb{Z}^{n}\}$, see \cite{FJ90}.

For any $\gamma \in \mathbb{Z}$, we put%
\begin{equation*}
\left\Vert f\right\Vert _{\widetilde{B}_{p(\cdot ),q(\cdot )}^{\alpha (\cdot
),p(\cdot )}}^{\ast }:=\sup_{\{P\in \mathcal{Q},|P|\leq 1\}}\left\Vert
\left( \frac{2^{v\alpha \left( \cdot \right) }\varphi _{v}\ast f}{%
|P|^{1/p(\cdot )}}\chi _{P}\right) _{v\geq v_{P}-\gamma }\right\Vert _{\ell
^{q(\cdot )}(L^{p(\cdot )})}<\infty
\end{equation*}%
and 
\begin{equation*}
\left\Vert f\right\Vert _{B_{p(\cdot ),q(\cdot )}^{\alpha (\cdot ),\tau
(\cdot )}}^{\ast }:=\sup_{P\in \mathcal{Q}}\left\Vert \left( \frac{%
2^{v\alpha \left( \cdot \right) }\varphi _{v}\ast f}{\left\Vert \chi
_{P}\right\Vert _{\tau (\cdot )}}\chi _{P}\right) _{v\geq v_{P}^{+}-\gamma
}\right\Vert _{\ell ^{q(\cdot )}(L^{p(\cdot )})}<\infty ,
\end{equation*}%
where $\varphi _{-\gamma }$ is replaced by $\Phi _{-\gamma }$.

\begin{lemma}
\label{new-equinorm3}Let $\alpha \in C_{\mathrm{loc}}^{\log }$, $p,q,\mathbb{%
\tau }\in \mathcal{P}_{0}^{\log }$ and $0<q^{+}<\infty $. The quasi-norms $%
\left\Vert f\right\Vert _{A_{p(\cdot ),q(\cdot )}^{\alpha (\cdot ),\tau
(\cdot )}}^{\ast }$ and $\left\Vert f\right\Vert _{A_{p(\cdot ),q(\cdot
)}^{\alpha (\cdot ),\tau (\cdot )}}$ are equivalent with equivalent
constants depending on $\gamma $.
\end{lemma}

\emph{Proof.} By similarity, we only consider $B_{p(\cdot ),q(\cdot
)}^{\alpha (\cdot ),\tau (\cdot )}$ and the case $\gamma >0$. First let us
prove that $\left\Vert f\right\Vert _{B_{p(\cdot ),q(\cdot )}^{\alpha (\cdot
),\tau (\cdot )}}^{\ast }\leq c\left\Vert f\right\Vert _{B_{p(\cdot
),q(\cdot )}^{\alpha (\cdot ),\tau (\cdot )}}$. By the scaling argument, it
suffices to consider the case $\left\Vert f\right\Vert _{B_{p(\cdot
),q(\cdot )}^{\alpha (\cdot ),\tau (\cdot )}}=1$ and show that the modular
of $f$ on the left-hand side is bounded. In particular, we will show that 
\begin{equation*}
\sum_{v=v_{P}^{+}-\gamma }^{\infty }\left\Vert \left\vert \frac{2^{v\alpha
(\cdot )}\varphi _{v}\ast f}{\left\Vert \chi _{P}\right\Vert _{\tau (\cdot )}%
}\right\vert ^{q(\cdot )}\chi _{P}\right\Vert _{\frac{p(\cdot )}{q(\cdot )}%
}\leq c
\end{equation*}%
for any dyadic cube $P$. As in \cite[Lemma 2.6]{WYY}, it suffices to prove
that for all dyadic cube $P$ with $l(P)\geq 1$,%
\begin{equation*}
I_{P}=\sum_{v=-\gamma }^{0}\left\Vert \left\vert \frac{2^{v\alpha (\cdot
)}\varphi _{v}\ast f}{\left\Vert \chi _{P}\right\Vert _{\tau (\cdot )}}%
\right\vert ^{q(\cdot )}\chi _{P}\right\Vert _{\frac{p(\cdot )}{q(\cdot )}%
}\leq c
\end{equation*}%
and for all dyadic cube $P$ with $l(P)<1$,%
\begin{equation*}
J_{P}=\sum_{v=v_{P}-\gamma }^{v_{P}-1}\left\Vert \left\vert \frac{2^{v\alpha
(\cdot )}\varphi _{v}\ast f}{\left\Vert \chi _{P}\right\Vert _{\tau (\cdot )}%
}\right\vert ^{q(\cdot )}\chi _{P}\right\Vert _{\frac{p(\cdot )}{q(\cdot )}%
}\leq c.
\end{equation*}%
The estimate of $I_{P}$, clearly follows from the inequality $\left\Vert
\left\vert \frac{\varphi _{v}\ast f}{\left\Vert \chi _{P}\right\Vert _{\tau
(\cdot )}}\right\vert ^{q(\cdot )}\chi _{P}\right\Vert _{\frac{p(\cdot )}{%
q(\cdot )}}\leq c$ for any $v=-\gamma ,...,0$ and any dyadic cube $P$ with $%
l(P)\geq 1$. This claim can be reformulated as showing that%
\begin{equation}
\left\Vert \frac{\varphi _{v}\ast f}{\left\Vert \chi _{P}\right\Vert _{\tau
(\cdot )}}\chi _{P}\right\Vert _{p(\cdot )}\leq c.  \label{estJ1}
\end{equation}%
By $\mathrm{\eqref{Ass1}}$ and $\mathrm{\eqref{Ass2}}$, there exist $\omega
_{v}\in \mathcal{S}(\mathbb{R}^{n})$, $v=-\gamma ,\cdot \cdot \cdot ,-1$ and 
$\eta _{1},\eta _{2}\in \mathcal{S}(\mathbb{R}^{n})$ such that%
\begin{equation*}
\varphi _{v}=\omega _{v}\ast \Phi ,\quad v=-\gamma ,\cdot \cdot \cdot ,-1%
\text{\quad and\quad }\varphi =\varphi _{0}=\eta _{1}\ast \Phi +\eta
_{2}\ast \varphi _{1}.
\end{equation*}%
Hence $\varphi _{v}\ast f=\omega _{v}\ast \Phi \ast f$ for $v=-\gamma
,...,-1 $ and $\varphi _{0}\ast f=\eta _{1}\ast \Phi \ast f+\eta _{2}\ast
\varphi _{1}\ast f$. Applying Lemma \ref{key-estimate1}, $\mathrm{%
\eqref{mod-est}}$ and the fact that $\left\Vert f\right\Vert _{B_{p(\cdot
),q(\cdot )}^{\alpha (\cdot ),\tau (\cdot )}}\leq 1$ to estimate the
left-hand side of $\mathrm{\eqref{estJ1}}$ by 
\begin{equation*}
C\left\Vert \Phi \ast f\right\Vert _{L_{\tau (\cdot )}^{p(\cdot
)}}+C\left\Vert \varphi _{1}\ast f\right\Vert _{L_{\tau (\cdot )}^{p(\cdot
)}}\leq c.
\end{equation*}%
To estimate $J_{P}$, denote by $P(\gamma )$ the dyadic cube containing $P$
with $l(P(\gamma ))=2^{\gamma }l(P)$. If $v_{P}\geq \gamma +1$, applying the
fact that $v_{P(\gamma )}=v_{P}-\gamma $, $\frac{\left\Vert \chi _{P(\gamma
)}\right\Vert _{\tau (\cdot )}}{\left\Vert \chi _{P}\right\Vert _{\tau
(\cdot )}}\approx c$ (see Lemma \ref{n-cube-est}) and $P\subset P(\gamma )$,
we then have%
\begin{equation*}
J_{P}\lesssim \sum_{v=v_{P(\gamma )}}^{v_{P}-1}\left\Vert \left\vert \frac{%
2^{v\alpha (\cdot )}\varphi _{v}\ast f}{\left\Vert \chi _{P(\gamma
)}\right\Vert _{\tau (\cdot )}}\right\vert ^{q(\cdot )}\chi _{P(\gamma
)}\right\Vert _{\frac{p(\cdot )}{q(\cdot )}}\leq c.
\end{equation*}%
If $1\leq v_{P}\leq \gamma $, we write $J_{P}=\sum_{v=v_{P}-\gamma
}^{-1}\cdot \cdot \cdot +\sum_{v=0}^{v_{P}-1}\cdot \cdot \cdot
=J_{P}^{1}+J_{P}^{2}$. Let $P(2^{v_{P}})$ the dyadic cube containing $P$
with $l(P(2^{v_{P}}))=2^{v_{P}}l(P)=1$, by the fact that $\frac{\left\Vert
\chi _{P(2^{v_{P}})}\right\Vert _{\tau (\cdot )}}{\left\Vert \chi
_{P}\right\Vert _{\tau (\cdot )}}\lesssim 2^{nv_{P}/\tau ^{-}}\leq
2^{n\gamma /\tau ^{-}}$, see Lemma \ref{n-cube-est}, we have%
\begin{equation*}
J_{P}^{2}\lesssim \sum_{v=v_{P(2^{v_{P}})}}^{v_{P}-1}\left\Vert \left\vert 
\frac{2^{v\alpha (\cdot )}\varphi _{v}\ast f}{\left\Vert \chi
_{P(2^{v_{P}})}\right\Vert _{\tau (\cdot )}}\right\vert ^{q(\cdot )}\chi
_{P(2^{v_{P}})}\right\Vert _{\frac{p(\cdot )}{q(\cdot )}}\leq c.
\end{equation*}%
By a similar argument to the estimate for $I_{P}$, we see that $%
J_{P}^{1}\leq c$.

For the converse estimate, it suffices to show that 
\begin{equation*}
\left\Vert \left\vert \frac{\Phi \ast f}{\left\Vert \chi _{P}\right\Vert
_{\tau (\cdot )}}\right\vert ^{q(\cdot )}\chi _{P}\right\Vert _{\frac{%
p(\cdot )}{q(\cdot )}}\leq c
\end{equation*}%
for all $P\in \mathcal{Q}$ with $l(P)\geq 1$ and all $f\in B_{p(\cdot
),q(\cdot )}^{\alpha (\cdot ),\tau (\cdot )}$ with $\left\Vert f\right\Vert
_{B_{p(\cdot ),q(\cdot )}^{\alpha (\cdot ),\tau (\cdot )}}^{\ast }\leq 1$.
This claim can be reformulated as showing that $\left\Vert \frac{\Phi \ast f%
}{\left\Vert \chi _{P}\right\Vert _{\tau (\cdot )}}\chi _{P}\right\Vert
_{p(\cdot )}\leq c$. Using the fact that there exist $\rho _{v}\in \mathcal{S%
}(\mathbb{R}^{n})$, $v=-\gamma ,$\textperiodcentered \textperiodcentered
\textperiodcentered $,1$, such that $\Phi \ast f=\rho _{-\gamma }\ast \Phi
_{-\gamma }\ast f+\sum_{v=1-\gamma }^{1}\rho _{v}\ast \varphi _{v}\ast f$,
see \cite[p. 130]{FJ90}. Applying Lemma \ref{key-estimate1} we obtain%
\begin{equation*}
\left\Vert \rho _{-\gamma }\ast \Phi _{-\gamma }\ast f\right\Vert _{L_{\tau
(\cdot )}^{p(\cdot )}}\lesssim \left\Vert \Phi _{-\gamma }\ast f\right\Vert
_{L_{\tau (\cdot )}^{p(\cdot )}}\leq c,
\end{equation*}%
and%
\begin{equation*}
\left\Vert \rho _{v}\ast \varphi _{v}\ast f\right\Vert _{L_{\tau (\cdot
)}^{p(\cdot )}}\lesssim \left\Vert \varphi _{v}\ast f\right\Vert _{L_{\tau
(\cdot )}^{p(\cdot )}}\leq c,\quad v=1-\gamma ,\cdot \cdot \cdot ,1,
\end{equation*}%
by using $\mathrm{\eqref{mod-est}}$ and the fact that $\left\Vert
f\right\Vert _{B_{p(\cdot ),q(\cdot )}^{\alpha (\cdot ),\tau (\cdot
)}}^{\ast }\leq 1$. The proof is complete. \hspace*{\fill}\rule{3mm}{3mm}

\begin{definition}
\label{sequence-space}Let $p,q,\tau \in \mathcal{P}_{0}$ and let $\alpha :%
\mathbb{R}^{n}\rightarrow \mathbb{R}$. Then for all complex valued sequences 
$\lambda =\{\lambda _{v,m}\in \mathbb{C}:v\in \mathbb{N}_{0},m\in \mathbb{Z}%
^{n}\}$ we define%
\begin{equation*}
\widetilde{b}_{p\left( \cdot \right) ,q\left( \cdot \right) }^{\alpha \left(
\cdot \right) ,p(\cdot )}:=\Big\{\lambda :\left\Vert \lambda \right\Vert _{%
\widetilde{b}_{p\left( \cdot \right) ,q\left( \cdot \right) }^{\alpha \left(
\cdot \right) ,p(\cdot )}}<\infty \Big\},
\end{equation*}%
where%
\begin{equation*}
\left\Vert \lambda \right\Vert _{\widetilde{b}_{p\left( \cdot \right)
,q\left( \cdot \right) }^{\alpha \left( \cdot \right) ,p(\cdot
)}}:=\sup_{P\in \mathcal{Q}}\left\Vert \left( \frac{\sum\limits_{m\in 
\mathbb{Z}^{n}}2^{v(\alpha \left( \cdot \right) +n/2)}\lambda _{v,m}\chi
_{v,m}}{|P|^{1/p(\cdot )}}\chi _{P}\right) _{v\geq v_{P}^{+}}\right\Vert
_{\ell ^{q(\cdot )}(L^{p(\cdot )})}
\end{equation*}%
and 
\begin{equation*}
b_{p\left( \cdot \right) ,q\left( \cdot \right) }^{\alpha \left( \cdot
\right) ,\tau (\cdot )}:=\Big\{\lambda :\left\Vert \lambda \right\Vert
_{b_{p\left( \cdot \right) ,q\left( \cdot \right) }^{\alpha \left( \cdot
\right) ,\tau (\cdot )}}<\infty \Big\}
\end{equation*}%
where%
\begin{equation*}
\left\Vert \lambda \right\Vert _{b_{p\left( \cdot \right) ,q\left( \cdot
\right) }^{\alpha \left( \cdot \right) ,\tau (\cdot )}}:=\sup_{P\in \mathcal{%
Q}}\left\Vert \left( \frac{\sum\limits_{m\in \mathbb{Z}^{n}}2^{v(\alpha
\left( \cdot \right) +n/2)}\lambda _{v,m}\chi _{v,m}}{\left\Vert \chi
_{P}\right\Vert _{\tau (\cdot )}}\chi _{P}\right) _{v\geq
v_{P}^{+}}\right\Vert _{\ell ^{q(\cdot )}(L^{p(\cdot )})}.
\end{equation*}
\end{definition}

If we replace dyadic cubes $P$ by arbitrary balls $B_{J}$ of $\mathbb{R}^{n}$
with $J\in \mathbb{Z}$, we then obtain equivalent quasi-norms, where the
supremum is taken over all $J\in \mathbb{Z}$\ and all balls $B_{J}$\ of $%
\mathbb{R}^{n}$. In the definition of $\widetilde{b}_{p\left( \cdot \right)
,q\left( \cdot \right) }^{\alpha \left( \cdot \right) ,p(\cdot )}$ the
supremum can be taken over all dyadic cube $P$, with $|P|\leq 1$. Similarly,
we use $a_{p\left( \cdot \right) ,q\left( \cdot \right) }^{\alpha \left(
\cdot \right) ,\tau (\cdot )}$ to denote either $b_{p\left( \cdot \right)
,q\left( \cdot \right) }^{\alpha \left( \cdot \right) ,\tau (\cdot )}$ or $%
\widetilde{b}_{p\left( \cdot \right) ,q\left( \cdot \right) }^{\alpha \left(
\cdot \right) ,p(\cdot )}$. Let $\alpha _{0},\alpha _{1}\in C_{\mathrm{loc}%
}^{\log }$ and $p_{0},p_{1},q,\tau \in \mathcal{P}_{0}^{\log }$ with $%
0<q^{+}<\infty $. If ${\alpha }_{0}>{\alpha }_{1}$\ and ${\alpha }_{0}{(x)-}%
\frac{n}{p_{0}(x)}={\alpha }_{1}{(x)-}\frac{n}{p_{1}(x)}$ with $\Big(\frac{%
p_{0}}{p_{1}}\Big)^{-}<1$, then us in \cite{D5}, we can prove the following
Sobolev-type embeddings%
\begin{equation*}
a_{{p}_{0}{(\cdot )},q{(\cdot )}}^{{\alpha }_{0}{(\cdot ),}\tau (\cdot
)}\hookrightarrow a_{{p}_{1}{(\cdot )},q{(\cdot )}}^{{\alpha }_{1}{(\cdot ),}%
\tau (\cdot )}.
\end{equation*}

\begin{lemma}
Let $\alpha \in C_{\mathrm{loc}}^{\log }$, $p,q,\mathbb{\tau }\in \mathcal{P}%
_{0}^{\log }$, $0<q^{+}<\infty $, $v\in \mathbb{N}_{0},m\in \mathbb{Z}^{n}$, 
$x\in Q_{v,m}$ and $\lambda \in a_{p\left( \cdot \right) ,q\left( \cdot
\right) }^{\alpha \left( \cdot \right) ,\tau (\cdot )}$. Then there exists $%
c>0$ independent of $v$ and $m$ such that%
\begin{equation*}
|\lambda _{v,m}|\leq c\text{ }2^{-v(\alpha (x)+n/2)}\left\Vert \lambda
\right\Vert _{b_{p\left( \cdot \right) ,q\left( \cdot \right) }^{\alpha
\left( \cdot \right) ,\tau (\cdot )}}\left\Vert \chi _{v,m}\right\Vert
_{\tau (\cdot )}\left\Vert \chi _{v,m}\right\Vert _{p(\cdot )}^{-1}
\end{equation*}%
and 
\begin{equation*}
|\lambda _{v,m}|\leq c\text{ }2^{-v(\alpha (x)+n/2)}\left\Vert \lambda
\right\Vert _{\widetilde{b}_{p(\cdot ),q(\cdot )}^{\alpha (\cdot ),p(\cdot
)}}.
\end{equation*}
\end{lemma}

\emph{Proof.} By similarity, we only consider $\widetilde{b}_{p(\cdot
),q(\cdot )}^{\alpha (\cdot ),p(\cdot )}$. Let $\lambda \in \widetilde{b}%
_{p(\cdot ),q(\cdot )}^{\alpha (\cdot ),p(\cdot )},v\in \mathbb{N}_{0},m\in 
\mathbb{Z}^{n}$ and $x\in Q_{v,m}$, with $Q_{v,m}\in \mathcal{Q}$. Then $%
|\lambda _{v,m}|^{p^{-}}=|Q_{v,m}|^{-1}\int_{Q_{v,m}}|\lambda
_{v,m}|^{p^{-}}\chi _{v,m}(y)dy$. Using the fact that $2^{v(\alpha \left(
x\right) -\alpha \left( y\right) )}\leq c$ for any $x,y\in Q_{v,m}$ and $%
|Q_{v,m}|^{1/p(x)}\approx \left\Vert \chi _{v,m}\right\Vert _{p(\cdot )}$,
see $\mathrm{\eqref{DHHR1}}$, we obtain%
\begin{eqnarray*}
\frac{2^{v(\alpha \left( x\right) +n/2)p^{-}}}{|Q_{v,m}|^{p^{-}/p(x)}}%
|\lambda _{v,m}|^{p^{-}} &\lesssim &|Q_{v,m}|^{-1}\int_{Q_{v,m}}\frac{%
2^{v(\alpha \left( y\right) +n/2)p^{-}}}{\left\Vert \chi _{v,m}\right\Vert
_{p(\cdot )}^{p^{-}}}|\lambda _{v,m}|^{p^{-}}\chi _{v,m}(y)dy \\
&\lesssim &|Q_{v,m}|^{-1}\int_{Q_{v,m}}\frac{2^{v(\alpha \left( y\right)
+n/2)p^{-}}}{|Q_{v,m}|^{p^{-}/p(y)}}|\lambda _{v,m}|^{p^{-}}\chi _{v,m}(y)dy.
\end{eqnarray*}%
Applying H\"{o}lder's inequality to estimate this expression by 
\begin{eqnarray*}
&&c|Q_{v,m}|^{-1}\left\Vert \frac{2^{v(\alpha \left( \cdot \right)
+n/2)p^{-}}}{|Q_{v,m}|^{p^{-}/p(\cdot )}}|\lambda _{v,m}|^{p^{-}}\chi
_{v,m}\right\Vert _{p/p^{-}}\left\Vert \chi _{v,m}\right\Vert
_{(p/p^{-})^{\prime }} \\
&\lesssim &\left\Vert \lambda \right\Vert _{\widetilde{b}_{p\left( \cdot
\right) ,q\left( \cdot \right) }^{\alpha \left( \cdot \right) ,p(\cdot
)}}^{p^{-}}\left\Vert \chi _{v,m}\right\Vert _{p/p^{-}}^{-1},
\end{eqnarray*}%
where we have used $\mathrm{\eqref{DHHR}}$. Therefore for any $x\in Q_{v,m}$%
\begin{eqnarray*}
|\lambda _{v,m}| &\lesssim &\text{ }2^{-v(\alpha
(x)+n/2)}|Q_{v,m}|^{1/p(x)}\left\Vert \lambda \right\Vert _{\widetilde{b}%
_{p\left( \cdot \right) ,q\left( \cdot \right) }^{\alpha \left( \cdot
\right) ,p(\cdot )}}\left\Vert \chi _{v,m}\right\Vert _{p(\cdot )}^{-1} \\
&\lesssim &\text{ }2^{-v(\alpha (x)+n/2)}\left\Vert \lambda \right\Vert _{%
\widetilde{b}_{p\left( \cdot \right) ,q\left( \cdot \right) }^{\alpha \left(
\cdot \right) ,p(\cdot )}},
\end{eqnarray*}%
again by $\mathrm{\eqref{DHHR1}}$, which completes the proof. \hspace*{\fill}%
\rule{3mm}{3mm}

\begin{lemma}
Let $\alpha \in C_{\mathrm{loc}}^{\log }$, $p,q,\mathbb{\tau }\in \mathcal{P}%
_{0}^{\log }$ and $\Psi $, $\psi \in \mathcal{S}(\mathbb{R}^{n})$ satisfy,
respectively, $\mathrm{\eqref{Ass1}}$ and $\mathrm{\eqref{Ass2}}$. Then for
all $\lambda \in a_{p\left( \cdot \right) ,q\left( \cdot \right) }^{\alpha
\left( \cdot \right) ,\tau (\cdot )}$%
\begin{equation*}
T_{\psi }\lambda :=\sum_{m\in \mathbb{Z}^{n}}\lambda _{0,m}\Psi
_{m}+\sum_{v=1}^{\infty }\sum_{m\in \mathbb{Z}^{n}}\lambda _{v,m}\psi _{v,m},
\end{equation*}%
converges in $\mathcal{S}^{\prime }(\mathbb{R}^{n})$; moreover, $T_{\psi
}:a_{p\left( \cdot \right) ,q\left( \cdot \right) }^{\alpha \left( \cdot
\right) ,\tau (\cdot )}\rightarrow \mathcal{S}^{\prime }(\mathbb{R}^{n})$ is
continuous.
\end{lemma}

\emph{Proof.} By similarity, we only consider $b_{p(\cdot ),q(\cdot
)}^{\alpha (\cdot ),\tau (\cdot )}$. Let $\lambda \in b_{p\left( \cdot
\right) ,q\left( \cdot \right) }^{\alpha \left( \cdot \right) ,\tau (\cdot
)} $ and $\phi \in \mathcal{S}(\mathbb{R}^{n})$. Observe that 
\begin{eqnarray*}
|\lambda _{v,m}| &\lesssim &\frac{2^{-v(\alpha ^{-}+n/2)}\left\Vert \lambda
\right\Vert _{b_{p\left( \cdot \right) ,q\left( \cdot \right) }^{\alpha
\left( \cdot \right) ,\tau (\cdot )}}\left\Vert \chi _{v,m}\right\Vert
_{\tau (\cdot )}}{\left\Vert \chi _{v,m}\right\Vert _{p(\cdot )}} \\
&\lesssim &\text{ }2^{vn(1/p^{-}-1/2-\alpha ^{-}/n-1/\tau ^{+})}\left\Vert
\lambda \right\Vert _{b_{p\left( \cdot \right) ,q\left( \cdot \right)
}^{\alpha \left( \cdot \right) ,\tau (\cdot )}}
\end{eqnarray*}%
for all dyadic cubes $Q_{v,m}$. Let $M>\max (n,n/p^{-}-\alpha ^{-}-n/\tau
^{+}-n)$. We see that,%
\begin{eqnarray*}
\sum_{m\in \mathbb{Z}^{n}}|\lambda _{0,m}||\langle \Psi _{m},\phi \rangle |
&\lesssim &\left\Vert \lambda \right\Vert _{b_{p\left( \cdot \right)
,q\left( \cdot \right) }^{\alpha \left( \cdot \right) ,\tau (\cdot
)}}\sum_{m\in \mathbb{Z}^{n}}\int |\Psi (x-m)||\phi (x)|dx \\
&\lesssim &\left\Vert \lambda \right\Vert _{b_{p\left( \cdot \right)
,q\left( \cdot \right) }^{\alpha \left( \cdot \right) ,\tau (\cdot
)}}\left\Vert \Psi \right\Vert _{\mathcal{S}_{2M}}\left\Vert \phi
\right\Vert _{\mathcal{S}_{M}}\sum_{m\in \mathbb{Z}^{n}}(1+|m|)^{-n-M}.
\end{eqnarray*}%
On the other hand, by \cite[Lemma 2.4]{WYY}, we obtain%
\begin{eqnarray*}
&&\sum_{v=1}^{\infty }\sum_{m\in \mathbb{Z}^{n}}|\lambda _{v,m}||\langle
\psi _{v,m},\phi \rangle | \\
&\lesssim &\left\Vert \psi \right\Vert _{\mathcal{S}_{M+1}}\left\Vert \phi
\right\Vert _{\mathcal{S}_{M+1}}\left\Vert \lambda \right\Vert _{b_{p\left(
\cdot \right) ,q\left( \cdot \right) }^{\alpha \left( \cdot \right) ,\tau
(\cdot )}}\sum_{v=1}^{\infty }\sum_{m\in \mathbb{Z}^{n}}\frac{2^{-vn(\alpha
^{-}/n+1/\tau ^{+}-1/p^{-}+1+M/n)}}{(1+|2^{-v}m|)^{n+M}} \\
&\lesssim &\left\Vert \psi \right\Vert _{\mathcal{S}_{M+1}}\left\Vert \phi
\right\Vert _{\mathcal{S}_{M+1}}\left\Vert \lambda \right\Vert _{b_{p\left(
\cdot \right) ,q\left( \cdot \right) }^{\alpha \left( \cdot \right) ,\tau
(\cdot )}},
\end{eqnarray*}%
which completes the proof. \hspace*{\fill}\rule{3mm}{3mm}

For a sequence $\lambda =\{\lambda _{v,m}\in \mathbb{C}:v\in \mathbb{N}%
_{0},m\in \mathbb{Z}^{n}\},0<r\leq \infty $ and a fixed $d>0$, set%
\begin{equation*}
\lambda _{v,m,r,d}^{\ast }:=\left( \sum_{h\in \mathbb{Z}^{n}}\frac{|\lambda
_{v,h}|^{r}}{(1+2^{v}|2^{-v}h-2^{-v}m|)^{d}}\right) ^{1/r}
\end{equation*}%
and $\lambda _{r,d}^{\ast }:=\{\lambda _{v,m,r,d}^{\ast }\in \mathbb{C}:v\in 
\mathbb{N}_{0},m\in \mathbb{Z}^{n}\}$.

\begin{lemma}
\label{lamda-equi}Let $\alpha \in C_{\mathrm{loc}}^{\log }$, $p,q,\mathbb{%
\tau }\in \mathcal{P}_{0}^{\log }$, $0<q^{+}<\infty $, $0<r<p^{-}$ and $%
a=r\max (2c_{\log }(q)+c_{\log }(\alpha ),2(\frac{1}{q^{-}}-\frac{1}{q^{+}}%
)+\alpha ^{+}-\alpha ^{-})$. Then%
\begin{equation*}
\left\Vert \lambda _{r,d}^{\ast }\right\Vert _{a_{p\left( \cdot \right)
,q\left( \cdot \right) }^{\alpha \left( \cdot \right) ,\tau (\cdot
)}}\approx \left\Vert \lambda \right\Vert _{a_{p\left( \cdot \right)
,q\left( \cdot \right) }^{\alpha \left( \cdot \right) ,\tau (\cdot )}}
\end{equation*}%
where 
\begin{equation*}
d>\left\{ 
\begin{array}{ccc}
n+a+n/\tau ^{-} & \text{if} & a_{p\left( \cdot \right) ,q\left( \cdot
\right) }^{\alpha \left( \cdot \right) ,\tau (\cdot )}=b_{p\left( \cdot
\right) ,q\left( \cdot \right) }^{\alpha \left( \cdot \right) ,\tau (\cdot )}
\\ 
n+a+c_{\log }(1/p)+n/p^{-} & \text{if} & a_{p\left( \cdot \right) ,q\left(
\cdot \right) }^{\alpha \left( \cdot \right) ,\tau (\cdot )}=\tilde{b}%
_{p\left( \cdot \right) ,q\left( \cdot \right) }^{\alpha \left( \cdot
\right) ,p(\cdot )}.%
\end{array}%
\right.
\end{equation*}
\end{lemma}

The proof of this lemma is postponed to the Appendix.

\begin{theorem}
\label{phi-tran}Let $\alpha \in C_{\mathrm{loc}}^{\log }$, $p,q,\mathbb{\tau 
}\in \mathcal{P}_{0}^{\log }$ and $0<q^{+}<\infty $. \textit{Suppose that }$%
\Phi $, $\Psi \in \mathcal{S}(\mathbb{R}^{n})$ satisfying $\mathrm{%
\eqref{Ass1}}$ and $\varphi ,\psi \in \mathcal{S}(\mathbb{R}^{n})$ satisfy $%
\mathrm{\eqref{Ass2}}$ such that $\mathrm{\eqref{Ass4}}$ holds. The
operators $S_{\varphi }:A_{p\left( \cdot \right) ,q\left( \cdot \right)
}^{\alpha \left( \cdot \right) ,\tau (\cdot )}\rightarrow a_{p\left( \cdot
\right) ,q\left( \cdot \right) }^{\alpha \left( \cdot \right) ,\tau (\cdot
)} $ and $T_{\psi }:a_{p\left( \cdot \right) ,q\left( \cdot \right)
}^{\alpha \left( \cdot \right) ,\tau (\cdot )}\rightarrow A_{p\left( \cdot
\right) ,q\left( \cdot \right) }^{\alpha \left( \cdot \right) ,\tau (\cdot
)} $ are bounded. Furthermore, $T_{\psi }\circ S_{\varphi }$ is the identity
on $A_{p\left( \cdot \right) ,q\left( \cdot \right) }^{\alpha \left( \cdot
\right) ,\tau (\cdot )}$.
\end{theorem}

{\emph{Proof.} By similarity, we only consider $b_{p(\cdot ),q(\cdot
)}^{\alpha (\cdot ),\tau (\cdot )}$ and $B_{p(\cdot ),q(\cdot )}^{\alpha
(\cdot ),\tau (\cdot )}$. For any $f\in \mathcal{S}^{\prime }(\mathbb{R}%
^{n}) $ we put $\sup (f):=\{\sup_{v,m}(f):v\in \mathbb{N}_{0},m\in \mathbb{Z}%
^{n}\} $ where 
\begin{equation*}
\sup_{v,m}(f):=2^{-vn/2}\sup_{y\in Q_{v,m}}|\widetilde{\varphi _{v}}\ast
f(y)|
\end{equation*}%
if $v\in \mathbb{N},m\in \mathbb{Z}^{n}$ and 
\begin{equation*}
\sup_{0,m}(f):=\sup_{y\in Q_{0,m}}|\widetilde{\Phi }\ast f(y)|
\end{equation*}%
if $m\in \mathbb{Z}^{n}$. For any $\gamma >0$, we define the sequence $%
\inf_{\gamma }(f):=\{\inf_{v,m,\gamma }(f):v\in \mathbb{N}_{0},m\in \mathbb{Z%
}^{n}\}$ by setting 
\begin{equation*}
\inf_{v,m,\gamma }(f):=2^{-vn/2}\sup_{h\in \mathbb{Z}^{n}}\{\inf_{y\in
Q_{v+\gamma ,h}}|\widetilde{\varphi _{v}}\ast f(y)|:Q_{v+\gamma ,h}\cap
Q_{v,m}\neq \emptyset \}
\end{equation*}%
if $v\in \mathbb{N},m\in \mathbb{Z}^{n}$ and 
\begin{equation*}
\inf_{0,m,\gamma }(f):=\sup_{h\in \mathbb{Z}^{n}}\{\inf_{y\in Q_{\gamma ,h}}|%
\widetilde{\Phi }\ast f(y)|:Q_{\gamma ,h}\cap Q_{0,m}\neq \emptyset \}
\end{equation*}%
if $m\in \mathbb{Z}^{n}$. Here $\widetilde{\varphi _{j}}(x):=2^{jn}\overline{%
\varphi (-2^{j}x)}$ and $\widetilde{\Phi }(x):=\overline{\Phi (-x)}$. As in
Lemma A.5 of \cite{FJ90} we obtain 
\begin{equation*}
\left\Vert \text{inf}_{\gamma }(f)\right\Vert _{b_{p\left( \cdot \right)
,q\left( \cdot \right) }^{\alpha \left( \cdot \right) ,\tau (\cdot
)}}\lesssim \left\Vert f\right\Vert _{B_{p\left( \cdot \right) ,q\left(
\cdot \right) }^{\alpha \left( \cdot \right) ,\tau (\cdot )}}
\end{equation*}%
for any $\alpha \in C_{\text{loc}}^{\log }$, $p,q,\mathbb{\tau }\in \mathcal{%
P}_{0}^{\log }$, $0<q^{+}<\infty $ and $\gamma >0$ sufficiently large.
Indeed, we have%
\begin{equation*}
\left\Vert \text{inf}_{\gamma }(f)\right\Vert _{b_{p\left( \cdot \right)
,q\left( \cdot \right) }^{\alpha \left( \cdot \right) ,\tau (\cdot
)}}=c\sup_{P\in \mathcal{Q}}\left\Vert \left( \frac{\sum\limits_{m\in 
\mathbb{Z}^{n}}2^{j(\alpha \left( \cdot \right) +n/2)}\inf_{j-\gamma
,m,\gamma }(f)\chi _{j-\gamma ,m}}{\left\Vert \chi _{P}\right\Vert _{\tau
(\cdot )}}\chi _{P}\right) _{j\geq v_{P}^{+}+\gamma }\right\Vert _{\ell
^{q(\cdot )}(L^{p(\cdot )})}.
\end{equation*}%
Define a sequence $\{\lambda _{i,k}\}_{i\in \mathbb{N}_{0},k\in \mathbb{Z}%
^{n}}$ by setting $\lambda _{i,k}:=2^{-in/2}\inf_{y\in Q_{i,k}}|\widetilde{%
\varphi _{i-\gamma }}\ast f(y)|$ and $\lambda _{0,k}:=\inf_{y\in Q_{\gamma
,k}}|\widetilde{\Phi }\ast f(y)|$. We have%
\begin{equation*}
\inf_{j-\gamma ,m,\gamma }(f):=2^{\gamma n/2}\sup_{h\in \mathbb{Z}%
^{n}}\{\lambda _{j,h}:Q_{j,h}\cap Q_{j-\gamma ,m}\neq \emptyset \}
\end{equation*}%
and 
\begin{equation*}
\inf_{0,m,\gamma }(f):=\sup_{h\in \mathbb{Z}^{n}}\{\lambda _{0,h}:Q_{\gamma
,h}\cap Q_{0,m}\neq \emptyset \}.
\end{equation*}%
Let $h\in \mathbb{Z}^{n}$ with $Q_{j,h}\cap Q_{j-\gamma ,m}\neq \emptyset $.
Then $\lambda _{j,h}\leq c$ $2^{\gamma d/r}\lambda _{j,k,r,d}^{\ast }$ for
any $k\in \mathbb{Z}^{n}$ with $Q_{j,k}\cap Q_{j-\gamma ,m}\neq \emptyset $.
Hence%
\begin{equation*}
\sum\limits_{m\in \mathbb{Z}^{n}}\inf_{j-\gamma ,m,\gamma }(f)\chi
_{j-\gamma ,m}\lesssim \sum\limits_{k\in \mathbb{Z}^{n}}\lambda
_{j,k,r,d}^{\ast }\chi _{j,k}
\end{equation*}%
and%
\begin{equation*}
\left\Vert \text{inf}_{\gamma }(f)\right\Vert _{b_{p\left( \cdot \right)
,q\left( \cdot \right) }^{\alpha \left( \cdot \right) ,\tau (\cdot
)}}\lesssim \sup_{P\in \mathcal{Q}}\left\Vert \left( \frac{\sum\limits_{k\in 
\mathbb{Z}^{n}}2^{j(\alpha \left( \cdot \right) +n/2)}\lambda
_{j,k,r,d}^{\ast }\chi _{j,k}}{\left\Vert \chi _{P}\right\Vert _{\tau (\cdot
)}}\chi _{P}\right) _{j\geq v_{P}^{+}+\gamma }\right\Vert _{\ell ^{q(\cdot
)}(L^{p(\cdot )})}.
\end{equation*}%
Notice that $P$ $=\cup _{m=1}^{2^{\gamma n}}P_{m}$, where $%
\{P_{m}\}_{m=1}^{2^{\gamma n}}$ are disjoint dyadic cubes with side length $%
l(P_{m})=2^{-(v_{P}+\gamma )}$. Therefore, taking $0<s<\frac{1}{2}\min
(p^{-},q^{-},2)$ and applying Lemmas \ref{n-cube-est} and \ref{lamda-equi},%
\begin{eqnarray*}
\left\Vert \text{inf}_{\gamma }(f)\right\Vert _{b_{p\left( \cdot \right)
,q\left( \cdot \right) }^{\alpha \left( \cdot \right) ,\tau (\cdot )}}^{s}
&\lesssim &\sum_{m=1}^{2^{\gamma n}}\sup_{P\in \mathcal{Q}}\left\Vert \left( 
\frac{\sum\limits_{k\in \mathbb{Z}^{n}}2^{j(\alpha \left( \cdot \right)
+n/2)}\lambda _{j,k,r,d}^{\ast }\chi _{j,k}}{\left\Vert \chi _{P}\right\Vert
_{\tau (\cdot )}}\chi _{P_{m}}\right) _{j\geq (v_{P}+\gamma
)^{+}}\right\Vert _{\ell ^{q(\cdot )}(L^{p(\cdot )})}^{s} \\
&\lesssim &\sup_{P\in \mathcal{Q}}\left\Vert \left( \frac{\sum\limits_{k\in 
\mathbb{Z}^{n}}2^{j(\alpha \left( \cdot \right) +n/2)}\lambda _{j,k}\chi
_{j,k}}{\left\Vert \chi _{P}\right\Vert _{\tau (\cdot )}}\chi _{P}\right)
_{j\geq v_{P}^{+}}\right\Vert _{\ell ^{q(\cdot )}(L^{p(\cdot )})}^{s}.
\end{eqnarray*}%
By Lemma \ref{new-equinorm3}, we obtain%
\begin{equation*}
\left\Vert \text{inf}_{\gamma }(f)\right\Vert _{b_{p\left( \cdot \right)
,q\left( \cdot \right) }^{\alpha \left( \cdot \right) ,\tau (\cdot )}}\leq
c\left\Vert f\right\Vert _{B_{p(\cdot ),q(\cdot )}^{\alpha (\cdot ),\tau
(\cdot )}}^{\ast }\leq c\left\Vert f\right\Vert _{B_{p(\cdot ),q(\cdot
)}^{\alpha (\cdot ),\tau (\cdot )}}.
\end{equation*}%
Applying Lemma A.4 of \cite{FJ90}, see also Lemma 8.3 of \cite{M07}, we
obtain inf$_{\gamma }(f)_{r,d}^{\ast }\approx \sup (f)_{r,d}^{\ast }$. Hence
for $\gamma >0$ sufficiently large we obtain by applying Lemma \ref%
{lamda-equi}, $\left\Vert \text{inf}_{\gamma }(f)_{r,d}^{\ast }\right\Vert
_{b_{p\left( \cdot \right) ,q\left( \cdot \right) }^{\alpha \left( \cdot
\right) ,\tau (\cdot )}}\approx \left\Vert \text{inf}_{\gamma
}(f)\right\Vert _{b_{p\left( \cdot \right) ,q\left( \cdot \right) }^{\alpha
\left( \cdot \right) ,\tau (\cdot )}}$ and $\left\Vert \sup (f)_{r,d}^{\ast
}\right\Vert _{b_{p\left( \cdot \right) ,q\left( \cdot \right) }^{\alpha
\left( \cdot \right) ,\tau (\cdot )}}\approx \left\Vert \sup (f)\right\Vert
_{b_{p\left( \cdot \right) ,q\left( \cdot \right) }^{\alpha \left( \cdot
\right) ,\tau (\cdot )}}$ for any $\alpha \in C_{\text{loc}}^{\log }$, $p,q,%
\mathbb{\tau }\in \mathcal{P}_{0}^{\log }$, with $0<q^{+}<\infty $.
Therefore, 
\begin{equation*}
\left\Vert \text{inf}_{\gamma }(f)\right\Vert _{b_{p\left( \cdot \right)
,q\left( \cdot \right) }^{\alpha \left( \cdot \right) ,\tau (\cdot
)}}\approx \left\Vert f\right\Vert _{B_{p\left( \cdot \right) ,q\left( \cdot
\right) }^{\alpha \left( \cdot \right) ,\tau (\cdot )}}\approx \left\Vert
\sup (f)\right\Vert _{b_{p\left( \cdot \right) ,q\left( \cdot \right)
}^{\alpha \left( \cdot \right) ,\tau (\cdot )}}.
\end{equation*}%
Use these estimates and repeating the proof of Theorem 2.2 in \cite{FJ90} or
Theorem 2.1 in \cite{WYY}, and complete the proof of Theorem \ref{phi-tran}. 
\hspace*{\fill}\rule{3mm}{3mm} }

From Theorem \ref{phi-tran}, we obtain the next important property of spaces 
$A_{p\left( \cdot \right) ,q\left( \cdot \right) }^{\alpha \left( \cdot
\right) ,\tau (\cdot )}$.

\begin{corollary}
Let $\alpha \in C_{\mathrm{loc}}^{\log }$, $p,q,\mathbb{\tau }\in \mathcal{P}%
_{0}^{\log }$ and $0<q^{+}<\infty $, The definition of the spaces $%
A_{p\left( \cdot \right) ,q\left( \cdot \right) }^{\alpha \left( \cdot
\right) ,\tau (\cdot )}$ is independent of the choices of $\Phi $ and $%
\varphi $.
\end{corollary}

\section{Decomposition by atoms}

In recent years, it turned out that atomic and sub-atomic, as well as
wavelet decompositions of some function spaces are extremely useful in many
aspects. This concerns, for instance, the investigation of (compact)
embeddings between function spaces. But this applies equally to questions of
mapping properties of pseudo-differential operators and to trace problems,
where arguments can be equivalently transferred to the sequence space, which
is often more convenient to handle. The idea of atomic decompositions leads
back to M. Frazier and B. Jawerth in their series of papers \cite{FJ86}, 
\cite{FJ90}, see also \cite{T3}. \vskip5pt

The main goal of this section is to prove an atomic decomposition result for 
$B_{p(\cdot ),q(\cdot )}^{\alpha (\cdot ),\tau (\cdot )}$ and $\widetilde{B}%
_{p(\cdot ),q(\cdot )}^{\alpha (\cdot ),p(\cdot )}$. Atoms are the building
blocks for the atomic decomposition.

\begin{definition}
\label{Atom-Def}Let $K\in \mathbb{N}_{0},L+1\in \mathbb{N}_{0}$ and let $%
\gamma >1$. A $K$-times continuous differentiable function $a\in C^{K}(%
\mathbb{R}^{n})$ is called $[K,L]$-atom centered at $Q_{v,m}$, $v\in \mathbb{%
N}_{0}$ and $m\in \mathbb{Z}^{n}$, if

\begin{equation}
\mathrm{supp}\text{ }a\subseteq \gamma Q_{v,m}  \label{supp-cond}
\end{equation}

\begin{equation}
|\partial ^{\beta }a(x)|\leq 2^{v(|\beta |+1/2)}\text{,\quad for\quad }0\leq
|\beta |\leq K,x\in \mathbb{R}^{n}  \label{diff-cond}
\end{equation}%
and if%
\begin{equation}
\int_{\mathbb{R}^{n}}x^{\beta }a(x)dx=0,\text{\quad for\quad }0\leq |\beta
|\leq L\text{ and }v\geq 1.  \label{mom-cond}
\end{equation}
\end{definition}

If the atom $a$ located at $Q_{v,m}$, that means if it fulfills $\mathrm{%
\eqref{supp-cond}}$, then we will denote it by $a_{v,m}$. For $v=0$ or $L=-1$
there are no moment\ conditions\ $\mathrm{\eqref{mom-cond}}$ required.\vskip%
5pt

For proving the decomposition by atoms we need the following lemma, see
Frazier \& Jawerth \cite[Lemma 3.3]{FJ86}.

\begin{lemma}
\label{FJ-lemma}Let $\Phi $ and $\varphi $ satisfy, respectively, $\mathrm{%
\eqref{Ass1}}$ and $\mathrm{\eqref{Ass2}}$ and let $\rho _{v,m}$ be an $%
\left[ K,L\right] $-atom. Then%
\begin{equation*}
\left\vert \varphi _{j}\ast \rho _{v,m}(x)\right\vert \leq c\text{ }%
2^{(v-j)K+vn/2}\left( 1+2^{v}\left\vert x-x_{Q_{v,m}}\right\vert \right)
^{-M}
\end{equation*}%
if $v\leq j$, and%
\begin{equation*}
\left\vert \varphi _{j}\ast \rho _{v,m}(x)\right\vert \leq c\text{ }%
2^{(j-v)(L+n+1)+vn/2}\left( 1+2^{j}\left\vert x-x_{Q_{v,m}}\right\vert
\right) ^{-M}
\end{equation*}%
if $v\geq j$, where $M$ is sufficiently large, $\varphi _{j}=2^{jn}\varphi
(2^{j}\cdot )$ and $\varphi _{0}$ is replaced by $\Phi $.
\end{lemma}

Now we come to the atomic decomposition theorem.

\begin{theorem}
\label{atomic-dec}\textit{Let }$\alpha \in C_{\mathrm{loc}}^{\log }$\textit{%
\ and }$p,q,\tau \in \mathcal{P}_{0}^{\log }$ with $0<q^{-}\leq q^{+}<\infty 
$\textit{. }Let $0<p^{-}\leq p^{+}\leq \infty $ and let $K,L+1\in \mathbb{N}%
_{0}$ such that%
\begin{equation}
K\geq ([\alpha ^{+}+n/\tau ^{-}]+1)^{+},  \label{K,L,B-F-cond}
\end{equation}
respectively 
\begin{equation*}
K\geq ([\alpha ^{+}+n/p^{-}]+1)^{+}
\end{equation*}%
and%
\begin{equation}
L\geq \max (-1,[n(\frac{1}{\min (1,p^{-})}-1)-\alpha ^{-}]).
\label{K,L,B-cond}
\end{equation}%
Then $f\in \mathcal{S}^{\prime }(\mathbb{R}^{n})$ belongs to $B_{p(\cdot
),q(\cdot )}^{\alpha (\cdot ),\tau (\cdot )}$, respectively to $\widetilde{B}%
_{p(\cdot ),q(\cdot )}^{\alpha (\cdot ),p(\cdot )}$, if and only if it can
be represented as%
\begin{equation}
f=\sum\limits_{v=0}^{\infty }\sum\limits_{m\in \mathbb{Z}^{n}}\lambda
_{v,m}\varrho _{v,m},\text{ \ \ \ \ converging in }\mathcal{S}^{\prime }(%
\mathbb{R}^{n})\text{,}  \label{new-rep}
\end{equation}%
where $\varrho _{v,m}$ are $\left[ K,L\right] $-atoms and $\lambda
=\{\lambda _{v,m}\in \mathbb{C}:v\in \mathbb{N}_{0},m\in \mathbb{Z}^{n}\}\in
b_{p(\cdot ),q(\cdot )}^{\alpha (\cdot ),\tau (\cdot )}$, respectively $%
\lambda \in \widetilde{b}_{p(\cdot ),q(\cdot )}^{\alpha (\cdot ),p(\cdot )}$%
. Furthermore, $\mathrm{inf}\left\Vert \lambda \right\Vert _{b_{p(\cdot
),q(\cdot )}^{\alpha (\cdot ),\tau (\cdot )}}$, respectively $\mathrm{inf}%
\left\Vert \lambda \right\Vert _{\widetilde{b}_{p(\cdot ),q(\cdot )}^{\alpha
(\cdot ),\tau (\cdot )}}$, where the infimum is taken over admissible
representations\ $\mathrm{\eqref{new-rep}}$\textrm{, }is an equivalent
quasi-norm in $B_{p(\cdot ),q(\cdot )}^{\alpha (\cdot ),\tau (\cdot )}$,
respectively $\widetilde{B}_{p(\cdot ),q(\cdot )}^{\alpha (\cdot ),p(\cdot
)} $.
\end{theorem}

The convergence in $\mathcal{S}^{\prime }(\mathbb{R}^{n})$ can be obtained
as a by-product of the proof using the same method as in \cite[ Corollary
13.9 ]{T3}, so the convergence is postponed to the Appendix.\vskip5pt

If $p$, $q$, $\tau $, and $\alpha $ are constants, then the restriction $%
\mathrm{\eqref{K,L,B-F-cond}}$, and their counterparts, in the atomic
decomposition theorem are $K\geq ([\alpha +n/\tau ]+1)^{+}$\ and $L\geq \max
(-1,[n(\frac{1}{\min (1,p)}-1)-\alpha ])$, which are essentially the
restrictions from the works of \cite[Theorem 3.12]{D4}, with $\frac{1}{\tau }
$ in place of $\tau $.

{\emph{Proof.} By similarity, we only consider $B_{p(\cdot ),q(\cdot
)}^{\alpha (\cdot ),\tau (\cdot )}$. The proof follows the ideas in \cite[%
Theorem 6]{FJ86}.\vskip5pt }

\textit{Step }1\textit{.} Assume that $f\in B_{p\left( \cdot \right)
,q\left( \cdot \right) }^{\alpha \left( \cdot \right) ,\tau \left( \cdot
\right) }$. Using the same of the arguments used in \cite[Theorem 3]{D3} we
obtain a sequence $\{\lambda _{v,m}\}$ and $\rho _{v,m}$ (atoms in the sense
of Definition \ref{Atom-Def}) such that $f=\sum\limits_{v=0}^{\infty
}\sum\limits_{m\in \mathbb{Z}^{n}}\lambda _{v,m}\varrho _{v,m}$ and $%
\left\Vert \lambda \right\Vert _{b_{p\left( \cdot \right) ,q\left( \cdot
\right) }^{\alpha \left( \cdot \right) ,,\tau \left( \cdot \right) }}\leq
c\left\Vert f\right\Vert _{B_{p\left( \cdot \right) ,q\left( \cdot \right)
}^{\alpha \left( \cdot \right) ,\tau \left( \cdot \right) }}$.\vskip5pt

\textit{Step }2\textit{.} Assume that $f$ can be represented by $\mathrm{%
\eqref{new-rep}}$, with $K$ and $L$ satisfying $\mathrm{\eqref{K,L,B-F-cond}}
$ and $\mathrm{\eqref{K,L,B-cond}}$, respectively. We will show that $f\in
B_{p\left( \cdot \right) ,q\left( \cdot \right) }^{\alpha \left( \cdot
\right) ,\tau \left( \cdot \right) }$ and that for some $c>0$, $\left\Vert
f\right\Vert _{B_{p\left( \cdot \right) ,q\left( \cdot \right) }^{\alpha
\left( \cdot \right) ,\tau \left( \cdot \right) }}\leq c\left\Vert \lambda
\right\Vert _{b_{p\left( \cdot \right) ,q\left( \cdot \right) }^{\alpha
\left( \cdot \right) ,\tau \left( \cdot \right) }}$. We write 
\begin{equation*}
f=\sum\limits_{v=0}^{\infty }\sum\limits_{m\in \mathbb{Z}^{n}}\lambda
_{v,m}\rho _{v,m}=\sum\limits_{v=0}^{j}\cdot \cdot \cdot
+\sum\limits_{v=j+1}^{\infty }\cdot \cdot \cdot .
\end{equation*}%
Recalling the definition of $B_{p\left( \cdot \right) ,q\left( \cdot \right)
}^{\alpha \left( \cdot \right) ,\tau \left( \cdot \right) }$ space, it
suffices to estimate 
\begin{equation*}
\left( \sum\limits_{v=0}^{j}\sum\limits_{m\in \mathbb{Z}^{n}}2^{j\alpha
\left( \cdot \right) }\left\vert \lambda _{v,m}\right\vert \left\vert
\varphi _{j}\ast \rho _{v,m}\right\vert \right) _{j\geq 0}\text{ \ and \ }%
\left( \sum\limits_{v=j}^{\infty }\sum\limits_{m\in \mathbb{Z}%
^{n}}2^{j\alpha \left( \cdot \right) }\left\vert \lambda _{v,m}\right\vert
\left\vert \varphi _{j}\ast \rho _{v,m}\right\vert \right) _{j\geq 0}
\end{equation*}%
in $\ell ^{\tau (\cdot ),q(\cdot )}(L^{p(\cdot )})$-norm. From Lemma \ref%
{FJ-lemma}, we have for any $M$ sufficiently large and any $v\leq j$%
\begin{eqnarray*}
&&\sum\limits_{m\in \mathbb{Z}^{n}}2^{j\alpha \left( x\right) }\left\vert
\lambda _{v,m}\right\vert \left\vert \varphi _{j}\ast \rho
_{v,m}(x)\right\vert \\
&\lesssim &2^{(v-j)(K-\alpha ^{+})}\sum\limits_{m\in \mathbb{Z}%
^{n}}2^{v(\alpha \left( x\right) +n/2)}\left\vert \lambda _{v,m}\right\vert
\left( 1+2^{v}\left\vert x-x_{Q_{v,m}}\right\vert \right) ^{-M} \\
&=&2^{(v-j)(K-\alpha ^{+})}\sum\limits_{m\in \mathbb{Z}^{n}}2^{v(\alpha
\left( x\right) -n/2)}\left\vert \lambda _{v,m}\right\vert \eta
_{v,M}(x-x_{Q_{v,m}}) \\
&\lesssim &2^{(v-j)(K-\alpha ^{+})}\sum\limits_{m\in \mathbb{Z}%
^{n}}2^{v(\alpha \left( x\right) +1/2)}\left\vert \lambda _{v,m}\right\vert
\eta _{v,M}\ast \chi _{v,m}(x),
\end{eqnarray*}%
by Lemma \ref{Conv-est1}. Lemma \ref{DHR-lemma} gives $2^{v\alpha \left(
\cdot \right) }\eta _{v,M}\ast \chi _{v,m}\lesssim \eta _{v,T}\ast
2^{v\alpha \left( \cdot \right) }\chi _{v,m}$, with $T=M-c_{\log }(\alpha )$
and since $K>\alpha ^{+}+n/\tau ^{-}$ we apply Lemma \ref{Key-lemma} to
obtain%
\begin{eqnarray*}
&&\left\Vert \left( \sum_{v=0}^{j}2^{(v-j)(K-\alpha ^{+})}\eta _{v,T}\ast 
\left[ 2^{v(\alpha \left( \cdot \right) +n/2)}\sum\limits_{m\in \mathbb{Z}%
^{n}}\left\vert \lambda _{v,m}\right\vert \chi _{v,m}\right] \right)
_{j}\right\Vert _{\ell ^{\tau (\cdot ),q(\cdot )}(L^{p(\cdot )})} \\
&\lesssim &\left\Vert \left( \eta _{v,T}\ast \left[ 2^{v(\alpha \left( \cdot
\right) +n/2)}\sum\limits_{m\in \mathbb{Z}^{n}}\left\vert \lambda
_{v,m}\right\vert \chi _{v,m}\right] \right) _{v}\right\Vert _{\ell ^{\tau
(\cdot ),q(\cdot )}(L^{p(\cdot )})}.
\end{eqnarray*}%
The right-hand side can be rewritten us%
\begin{eqnarray*}
&&\sup_{P\in \mathcal{Q}}\left\Vert \left( \frac{\left( \eta _{v,T}\ast %
\left[ 2^{v(\alpha \left( \cdot \right) +n/2)}\sum\limits_{m\in \mathbb{Z}%
^{n}}\left\vert \lambda _{v,m}\right\vert \chi _{v,m}\right] \right) ^{r}}{%
\left\Vert \chi _{P}\right\Vert _{\tau (\cdot )}^{r}}\chi _{P}\right)
_{v\geq v_{P}^{+}}\right\Vert _{\ell ^{q(\cdot )/r}(L^{p(\cdot )/r})}^{1/r}
\\
&\lesssim &\sup_{P\in \mathcal{Q}}\left\Vert \left( \frac{\eta _{v,Tr}\ast %
\left[ 2^{v(\alpha \left( \cdot \right) +n/2)r}\sum\limits_{m\in \mathbb{Z}%
^{n}}\left\vert \lambda _{v,m}\right\vert ^{r}\chi _{v,m}\right] }{%
\left\Vert \chi _{P}\right\Vert _{\tau (\cdot )}^{r}}\chi _{P}\right)
_{v\geq v_{P}^{+}}\right\Vert _{\ell ^{q(\cdot )/r}(L^{p(\cdot )/r})}^{1/r},
\end{eqnarray*}%
by Lemma \ref{Conv-est2}, since $\eta _{v,T}\approx \eta _{v,T}\ast \eta
_{v,T}$ and $0<r<\min (1,p^{-})$. The application of Lemma \ref%
{Alm-Hastolemma1} and the fact that $\left\Vert \left( g_{v}\right) _{v\geq
v_{P}^{+}}\right\Vert _{\ell ^{q(\cdot )/r}(L^{p(\cdot
)/r})}^{1/r}=\left\Vert \left( |g_{v}|^{1/r}\right) _{v\geq
v_{P}^{+}}\right\Vert _{\ell ^{q(\cdot )}(L^{p(\cdot )})}$ give that the
last expression is bounded by $\left\Vert \lambda \right\Vert _{b_{p\left(
\cdot \right) ,q\left( \cdot \right) }^{\alpha \left( \cdot \right) ,\tau
\left( \cdot \right) }}$. Now from Lemma \ref{FJ-lemma}, we have for any $M$
sufficiently large and $v\geq j$%
\begin{eqnarray*}
&&\sum\limits_{m\in \mathbb{Z}^{n}}2^{j\alpha \left( x\right) }\left\vert
\lambda _{v,m}\right\vert \left\vert \varphi _{j}\ast \rho
_{v,m}(x)\right\vert \\
&\lesssim &2^{(j-v)(L+1+n/2)}\sum\limits_{m\in \mathbb{Z}^{n}}2^{j(\alpha
\left( x\right) +n/2)}\left\vert \lambda _{v,m}\right\vert \left(
1+2^{j}\left\vert x-x_{Q_{v,m}}\right\vert \right) ^{-M} \\
&=&2^{(j-v)(L+1+n/2)}\sum\limits_{m\in \mathbb{Z}^{n}}2^{j(\alpha \left(
x\right) -n/2)}\left\vert \lambda _{v,m}\right\vert \eta
_{j,M}(x-x_{Q_{v,m}}) \\
&\lesssim &2^{(j-v)(L+1+n/2)}\sum\limits_{m\in \mathbb{Z}^{n}}2^{j(\alpha
\left( x\right) -n/2)}\left\vert \lambda _{v,m}\right\vert \eta _{j,M}\ast
\eta _{v,M}(x-x_{Q_{v,m}}),
\end{eqnarray*}%
where the last inequality follows by Lemma \ref{Conv-est}, since $\eta
_{j,M}=\eta _{\min (v,j),M}$. Again by Lemma \ref{Conv-est1}, we have%
\begin{equation*}
\eta _{j,M}\ast \eta _{v,M}(x-x_{Q_{v,m}})\lesssim 2^{vn}\eta _{j,M}\ast
\eta _{v,M}\ast \chi _{v,m}(x).
\end{equation*}%
Therefore, $\sum\limits_{m\in \mathbb{Z}^{n}}2^{j\alpha \left( x\right)
}\left\vert \lambda _{v,m}\right\vert \left\vert \varphi _{j}\ast \rho
_{v,m}(x)\right\vert $ is bounded by%
\begin{eqnarray*}
&&c\text{ }2^{(j-v)(L+1-n/2)}\sum\limits_{m\in \mathbb{Z}^{n}}2^{j(\alpha
\left( x\right) +n/2)}\left\vert \lambda _{v,m}\right\vert \eta _{j,M}\ast
\eta _{v,M}\ast \chi _{v,m}(x) \\
&\lesssim &2^{(j-v)(L+1-\alpha ^{-})}\eta _{j,T}\ast \eta _{v,T}\ast \left[
2^{v(\alpha \left( \cdot \right) +n/2)}\sum\limits_{m\in \mathbb{Z}%
^{n}}\left\vert \lambda _{v,m}\right\vert \chi _{v,m}\right] (x),
\end{eqnarray*}%
by Lemma \ref{DHR-lemma}, with $T=M-c_{\log }(\alpha )$. Let $0<r<\min
(1,p^{-})$\ be a real number such that $L>n/r-1-\alpha ^{-}-n$. We have%
\begin{eqnarray*}
&&\left( \sum_{v=j}^{\infty }2^{(j-v)(L+1-\alpha ^{-})}\eta _{j,T}\ast \eta
_{v,T}\ast \left[ 2^{v(\alpha \left( \cdot \right) +n/2)}\sum\limits_{m\in 
\mathbb{Z}^{n}}\left\vert \lambda _{v,m}\right\vert \chi _{v,m}\right]
\right) ^{r} \\
&\leq &\sum_{v=j}^{\infty }2^{(j-v)(L+1-\alpha ^{-})r}\left( \eta _{j,T}\ast
\eta _{v,T}\ast \left[ 2^{v(\alpha \left( \cdot \right)
+n/2)}\sum\limits_{m\in \mathbb{Z}^{n}}\left\vert \lambda _{v,m}\right\vert
\chi _{v,m}\right] \right) ^{r} \\
&\leq &\sum_{v=j}^{\infty }2^{(j-v)(L-n/r+1-\alpha ^{-}+n)r}\eta _{j,Tr}\ast
\eta _{v,Tr}\ast \left[ 2^{v(\alpha \left( \cdot \right)
+n/2)r}\sum\limits_{m\in \mathbb{Z}^{n}}\left\vert \lambda _{v,m}\right\vert
^{r}\chi _{v,m}\right] ,
\end{eqnarray*}%
where the first estimate follows by the well-known inequality $\left(
\sum_{j=0}^{\infty }\left\vert a_{j}\right\vert \right) ^{\sigma }\leq
\sum_{j=0}^{\infty }\left\vert a_{j}\right\vert ^{\sigma }$, with $\left\{
a_{j}\right\} _{j}\subset \mathbb{C}$, $\sigma \in \left[ 0,1\right] $ and
the second inequality is by Lemma \ref{Conv-est2}. The application of Lemma %
\ref{Alm-Hastolemma1} gives that%
\begin{equation*}
\left\Vert \left( \sum_{v=j}^{\infty }2^{(j-v)(L+1-\alpha ^{-})}\eta
_{j,T}\ast \eta _{v,T}\ast \left[ 2^{v(\alpha \left( \cdot \right)
+n/2)}\sum\limits_{m\in \mathbb{Z}^{n}}\left\vert \lambda _{v,m}\right\vert
\chi _{v,m}\right] \right) _{j}\right\Vert _{\ell ^{\tau (\cdot ),q(\cdot
)}(L^{p(\cdot )})}
\end{equation*}%
is bounded by%
\begin{equation*}
c\sup_{P\in \mathcal{Q}}\left\Vert \left( \frac{\sum_{v=j}^{\infty
}2^{(j-v)Hr}\eta _{v,Tr}\ast \left[ 2^{v(\alpha \left( \cdot \right)
+n/2)r}\sum\limits_{m\in \mathbb{Z}^{n}}\left\vert \lambda _{vm}\right\vert
^{r}\chi _{v,m}\right] }{\left\Vert \chi _{P}\right\Vert _{\tau (\cdot )}^{r}%
}\chi _{P}\right) _{j\geq j_{P}^{+}}\right\Vert _{\ell ^{q(\cdot
)/r}(L^{p(\cdot )/r})}^{1/r},
\end{equation*}%
where $H:=L-n/r+n+1-\alpha ^{-}$. Observing that $H>0$, an application of
Lemma \ref{Key-lemma} (this is possible, see the proof of this lemma) yields
that the last expression is bounded by%
\begin{equation*}
c\sup_{P\in \mathcal{Q}}\left\Vert \left( \frac{\eta _{v,Tr}\ast \left[
2^{v(\alpha \left( \cdot \right) +n/2)r}\sum\limits_{m\in \mathbb{Z}%
^{n}}\left\vert \lambda _{v,m}\right\vert ^{r}\chi _{v,m}\right] }{%
\left\Vert \chi _{P}\right\Vert _{\tau (\cdot )}^{r}}\chi _{P}\right)
_{v\geq v_{P}^{+}}\right\Vert _{\ell ^{q(\cdot )/r}(L^{p(\cdot
)/r})}^{1/r}\lesssim \left\Vert \lambda \right\Vert _{b_{p\left( \cdot
\right) ,q\left( \cdot \right) }^{\alpha \left( \cdot \right) ,\tau \left(
\cdot \right) }},
\end{equation*}%
where we used again Lemma \ref{Alm-Hastolemma1} and hence the proof is
complete. \hspace*{\fill}\rule{3mm}{3mm}

\section{Appendix}

Here we present more technical proofs of the Lemmas.\vskip5pt

{\emph{Proof of Lemma} \ref{r-trick}. First let us prove that%
\begin{equation}
\left\vert \omega _{N}\ast g\left( x\right) \right\vert \leq c\text{ }(\eta
_{N,m}\ast \left\vert \omega _{N}\ast g\right\vert ^{r}(x))^{1/r},\quad x\in 
\mathbb{R}^{n},  \label{r-est}
\end{equation}%
where $c>0$ independent of $g$, $N$ and $x$. Let $\phi $ be a function in $%
\mathcal{S}\left( \mathbb{R}^{n}\right) $ satisfying $\mathcal{F}\phi =1$ on 
$\overline{B(0,1)}$. Then $\omega _{N}\ast g=\phi _{N}\ast \omega _{N}\ast g$
and we can distinguish two cases as follows: }

$\bullet $ $1\leq r<\infty $: observe that%
\begin{equation*}
\left\vert \omega _{N}\ast g\left( x\right) \right\vert \leq c\text{ }\int_{%
\mathbb{R}
^{n}}\eta _{N,m/r}(x-y)(1+N|x-y|)^{-m/r^{\prime }}\left\vert \omega _{N}\ast
g(y)\right\vert dy,
\end{equation*}%
where $r^{\prime }$ is the conjugate exponent of $r$ and we have used the
fact that $|\phi (Nx)|\leq c(1+N|x|)^{-m}$. By H\"{o}lder's inequality,%
\begin{eqnarray*}
\left\vert \omega _{N}\ast g\left( x\right) \right\vert &\leq &c\text{ }%
N^{n/r^{\prime }}(\eta _{N,m}\ast \left\vert \omega _{N}\ast g\right\vert
^{r}(x))^{1/r}\left\Vert (1+N|\cdot |)^{-m/r^{\prime }}\right\Vert
_{r^{\prime }} \\
&\leq &c\text{ }(\eta _{N,m}\ast \left\vert \omega _{N}\ast g\right\vert
^{r}(x))^{1/r}.
\end{eqnarray*}

$\bullet $ $0<r<1$: we put $g_{\omega ,N,m}^{\ast }(x):=\underset{y\in 
\mathbb{R}^{n}}{\sup }\frac{\left\vert \omega _{N}\ast g(y)\right\vert }{%
(1+N\left\vert y-x\right\vert )^{m}}$ and we have%
\begin{equation*}
\left\vert \omega _{N}\ast g(z)\right\vert \leq c\text{ }\int_{%
\mathbb{R}
^{n}}\eta _{N,m}(z-y)\left\vert \omega _{N}\ast g(y)\right\vert dy.
\end{equation*}%
We use the estimate $\left( 1+N\left\vert z-y\right\vert \right) ^{-m}\leq
\left( 1+N\left\vert z-x\right\vert \right) ^{m}\left( 1+N\left\vert
x-y\right\vert \right) ^{-m}$ we obtain%
\begin{eqnarray}
g_{\omega ,N,m}^{\ast }(x) &\leq &c\int_{%
\mathbb{R}
^{n}}\eta _{N,m}(x-y)\left\vert \omega _{N}\ast g(y)\right\vert
^{1-r}\left\vert \omega _{N}\ast g(y)\right\vert ^{r}dy  \notag \\
&\leq &c(g_{\omega ,N,m}^{\ast }(x))^{1-r}\int_{%
\mathbb{R}
^{n}}\eta _{N,m}(x-y)\left\vert \omega _{N}\ast g(y)\right\vert ^{r}dy\text{.%
}  \label{key-est1}
\end{eqnarray}%
Since $g$ is a tembered distribution and $\omega \in \mathcal{S}\left( 
\mathbb{R}^{n}\right) $, $|\omega _{N}\ast g(y)|$ is dominated by%
\begin{eqnarray*}
c\left\Vert \omega _{N}(y-\cdot )\right\Vert _{\mathcal{S}_{M}} &=&c\text{ }%
N^{n}\sup_{\gamma \in \mathbb{N}_{0}^{n},|\gamma |\leq M}\sup_{t\in \mathbb{R%
}^{n}}|\partial ^{\gamma }\omega (N(t-y))|(1+|t|)^{n+M+|\gamma |} \\
&\leq &c\text{ }\max (N^{M+n},N^{-2M})(1+N|y|)^{n+2M}
\end{eqnarray*}%
for some $M\in \mathbb{N}$ and any $y\in \mathbb{R}^{n}$, with $C>0$
independent of $N$ and $y$. Therefore,%
\begin{eqnarray*}
g_{\omega ,N,m}^{\ast }(x) &=&\underset{y\in \mathbb{R}^{n}}{\sup }\frac{%
\left\vert \omega _{N}\ast g(y)\right\vert }{(1+N\left\vert y-x\right\vert
)^{m}} \\
&\leq &\max (N^{M+n},N^{-2M})(1+N\left\vert y-x\right\vert
)^{-m}(1+N|y|)^{n+2M} \\
&\leq &C(N)(1+N\left\vert x\right\vert )^{n+2M}
\end{eqnarray*}%
if $m\geq n+2M$ and hence $g_{\omega ,N,m}^{\ast }(x)$ is finite (of course $%
m\geq n+2M$). Now we use the idea of \cite[Lemma 2.9]{Ry01}. Observe that
the right-hand side of $\mathrm{\eqref{r-est}}$ decreases as $m$ increases.
Therefore, we have $\mathrm{\eqref{r-est}}$ for all $m>n$ but with $c=c(g)$
depending on $g$. We can easily check that $\mathrm{\eqref{r-est}}$, with $%
c=c(g)$ imply that $g_{\omega ,N,m}^{\ast }(x)<\infty $. We assume that the
right-hand side of $\mathrm{\eqref{r-est}}$ is finite (otherwise, there is
nothing to prove). Returning to $\mathrm{\eqref{key-est1}}$ and having in
mind that now $g_{\omega ,N,m}^{\ast }(x)<\infty $, we end up with%
\begin{equation}
(g_{\omega ,N,m}^{\ast }(x))^{r}\leq c\int_{%
\mathbb{R}
^{n}}\eta _{N,m}(x-y)\left\vert \omega _{N}\ast g(y)\right\vert ^{r}dy\text{,%
}  \label{r-est1}
\end{equation}%
for all $m>n$ and $c$ independent of $g$, $N$ and $x$, which completes the
proof of $\mathrm{\eqref{r-est}}$. Now observe that%
\begin{eqnarray*}
\left\vert \theta _{R}\ast \omega _{N}\ast g\left( x\right) \right\vert 
&\leq &c\text{ }\int_{%
\mathbb{R}
^{n}}\eta _{R,d}(x-y)\left\vert \omega _{N}\ast g(y)\right\vert dy \\
&\leq &c\text{ }g_{\omega ,N,m}^{\ast }(x)\int_{%
\mathbb{R}
^{n}}\eta _{R,d}(x-y)(1+N|x-y|)^{m}dy \\
&\leq &c\text{ }\max \Big(1,\Big(\frac{N}{R}\Big)^{m}\Big)g_{\omega
,N,m}^{\ast }(x)R^{n}\int_{%
\mathbb{R}
^{n}}(1+R|x-y|)^{m-d}dy \\
&\leq &c\text{ }\max \Big(1,\Big(\frac{N}{R}\Big)^{m}\Big)g_{\omega
,N,m}^{\ast }(x),
\end{eqnarray*}%
provided we pick $d>m+n$ and $c$ independent of $g$, $N$ and $x$. Hence the
proof of $\mathrm{\eqref{r-trick-est}}$ is complete by using $\mathrm{%
\eqref{r-est1}}$.

{\emph{Proof of Lemma }\ref{key-estimate1}. By similarity, we only consider $%
L_{\tau (\cdot )}^{p(\cdot )}$. We use Lemma \ref{r-trick}, in the form%
\begin{equation*}
\left\vert \theta _{R}\ast \omega _{N}\ast f\left( x\right) \right\vert \leq
c\text{ }\max \Big(1,\Big(\frac{N}{R}\Big)^{m}\Big)(\eta _{N,m}\ast
\left\vert \omega _{N}\ast f\right\vert ^{r}(x))^{1/r}.
\end{equation*}%
where $0<r<p^{-}$, $m>2n+c_{\log }(\frac{1}{\tau })r$ and $x\in P$. We have,
with $k=(k_{1},...,k_{n})$, 
\begin{eqnarray*}
\eta _{N,m}\ast \left\vert \omega _{N}\ast f\right\vert ^{r}(x)
&=&N^{n}\int_{\mathbb{R}^{n}}\frac{|\omega _{N}\ast f(z)|^{r}}{\left(
1+N\left\vert x-z\right\vert \right) ^{m}}dz \\
&=&\int_{3P}\cdot \cdot \cdot dz+\sum_{k\in \mathbb{Z}^{n},%
\max_{i=1,...,n}|k_{i}|\geq 2}\int_{P+kl(P)}\cdot \cdot \cdot dz \\
&=&J_{N}^{1}(\omega _{N}\ast f)(x)+\sum_{k\in \mathbb{Z}^{n},%
\max_{i=1,...,n}|k_{i}|\geq 2}J_{N,k}^{2}(\omega _{N}\ast f)(x).
\end{eqnarray*}%
Thus we obtain%
\begin{eqnarray}
&&\left\Vert \frac{\theta _{R}\ast \omega _{N}\ast f}{\left\Vert \chi
_{P}\right\Vert _{\tau (\cdot )}}\chi _{P}\right\Vert _{p(\cdot )}^{r} 
\notag \\
&\lesssim &\max \Big(1,\Big(\frac{N}{R}\Big)^{mr}\Big)\left\Vert \frac{%
J_{N}^{1}(\omega _{N}\ast f)}{\left\Vert \chi _{P}\right\Vert _{\tau (\cdot
)}^{r}}\chi _{P}\right\Vert _{p(\cdot )/r}  \notag \\
&&+\max \Big(1,\Big(\frac{N}{R}\Big)^{mr}\Big)\sum_{k\in \mathbb{Z}%
^{n},\max_{i=1,...,n}|k_{i}|\geq 2}\left\Vert \frac{J_{N,k}^{2}(\omega
_{N}\ast f)}{\left\Vert \chi _{P}\right\Vert _{\tau (\cdot )}^{r}}\chi
_{P}\right\Vert _{p(\cdot )/r}.  \label{norm1}
\end{eqnarray}%
Let us prove that the first norm on the right-hand side is bounded by 
\begin{equation}
c\left\Vert \frac{\omega _{N}\ast f}{\left\Vert \chi _{P}\right\Vert _{\tau
(\cdot )}}\chi _{3P}\right\Vert _{p(\cdot )}^{r}.  \label{exp1}
\end{equation}%
We have%
\begin{equation*}
|J_{N}^{1}(\omega _{N}\ast f)(x)|\lesssim \text{ }N^{n}\int_{\mathbb{R}^{n}}%
\frac{\left\vert \omega _{N}\ast f(z)\right\vert ^{r}\chi _{3P}(z)}{\left(
1+N\left\vert x-z\right\vert \right) ^{m}}dz.
\end{equation*}%
Now the function $z\mapsto \frac{1}{\left( 1+\left\vert z\right\vert \right)
^{m}}$ is in $L^{1}$ (since $m>n$), then using the majorant property for the
Hardy-Littlewood maximal operator $\mathcal{M}$, see E. M. Stein and G.
Weiss \cite[Chapiter 2, (3.9)]{SW}, $\Big(|g|\ast \frac{1}{\left(
1+\left\vert \cdot \right\vert \right) ^{m}}\Big)(x)\lesssim $ $\left\Vert 
\frac{1}{\left( 1+\left\vert \cdot \right\vert \right) ^{m}}\right\Vert _{1}%
\mathcal{M}(g)(x)$, it follows that for any $x\in P$, $|J_{N}^{1}(\omega
_{N}\ast f)(x)|\leq C\text{ }\mathcal{M}\Big(|\omega _{N}\ast f|^{r}\chi
_{3P}\Big)(x)$ where the constant $C>0$ is independent of $x$ and $N$. Hence
the first norm of $\mathrm{\eqref{norm1}}$ is bounded by%
\begin{equation*}
c\left\Vert \mathcal{M}\Big(\frac{|\omega _{N}\ast f|^{r}}{\left\Vert \chi
_{P}\right\Vert _{\tau (\cdot )}^{r}}\chi _{3P}\Big)\right\Vert _{p(\cdot
)/r}\lesssim \left\Vert \frac{|\omega _{N}\ast f|^{r}}{\left\Vert \chi
_{P}\right\Vert _{\tau (\cdot )}^{r}}\chi _{3P}\right\Vert _{p(\cdot
)/r}=\left\Vert \frac{\omega _{N}\ast f}{\left\Vert \chi _{P}\right\Vert
_{\tau (\cdot )}}\chi _{3P}\right\Vert _{p(\cdot )}^{r},
\end{equation*}%
after using the fact that $\mathcal{M}:L^{p(\cdot )/r}\rightarrow L^{p(\cdot
)/r}$ is bounded. Notice that $3P=\cup _{h=1}^{3^{n}}P_{h}$, where $%
\{P_{h}\}_{h=1}^{3^{n}}$ are disjoint dyadic cubes with side length $%
l(P_{h})=l(P)$. Therefore $\chi _{3P}=\sum_{h=1}^{3^{n}}\chi _{P_{h}}$ and
the expression in $\mathrm{\eqref{exp1}}$ can be estimated by%
\begin{equation*}
c\sum_{h=1}^{3^{n}}\left\Vert \frac{\omega _{N}\ast f}{\left\Vert \chi
_{P_{h}}\right\Vert _{\tau (\cdot )}}\chi _{P_{h}}\right\Vert _{p(\cdot
)}^{r}\lesssim \left\Vert \omega _{N}\ast f\right\Vert _{L_{\tau (\cdot
)}^{p(\cdot )}}^{r},
\end{equation*}%
where we have used the fact that $\frac{\left\Vert \chi _{P_{h}}\right\Vert
_{\tau (\cdot )}}{\left\Vert \chi _{P}\right\Vert _{\tau (\cdot )}}\leq c$,
see Lemma \ref{n-cube-est} (ii) and the proof of the first part is finished.
The summation in $\mathrm{\eqref{norm1}}$ can be rewritten us%
\begin{equation}
\sum_{k\in \mathbb{Z}^{n},|k|\leq 4\sqrt{n}}\cdot \cdot \cdot +\sum_{k\in 
\mathbb{Z}^{n},|k|>4\sqrt{n}}\cdot \cdot \cdot .  \label{sec-sum}
\end{equation}%
The estimate of the first sum follows in the same manner as in the estimate
of $J_{N}^{1}(\omega _{N}\ast f)$, so we need only to estimate the second
sum. Let us prove that%
\begin{equation}
\left\Vert |k|^{m-n-c_{\log }(\frac{1}{\tau })r}\frac{J_{N,k}^{2}(\omega
_{N}\ast f)}{\left\Vert \chi _{P}\right\Vert _{\tau (\cdot )}^{r}}\chi
_{P}\right\Vert _{p(\cdot )/r}\lesssim \left( Nl(P)\right) ^{n-m}\left\Vert 
\frac{|\omega _{N}\ast f|^{r}}{\left\Vert \chi _{P+kl(P)}\right\Vert _{\tau
(\cdot )}^{r}}\chi _{P+kl(P)}\right\Vert _{p(\cdot )/r}.  \label{exp3}
\end{equation}%
Let $x\in P$, $z\in P+kl(P)$ with $k\in \mathbb{Z}^{n}$ and $|k|>4\sqrt{n}$.
Then $\left\vert x-z\right\vert \geq \frac{3}{4}\left\vert k\right\vert l(P)$
and the term $|J_{N,k}^{2}(\omega _{N}\ast f)(x)|$ is bounded by 
\begin{eqnarray*}
&&C\text{ }|k|^{-m}N^{n-m}\left( l(P)\right) ^{-m}\int_{P+kl(P)}\left\vert
\omega _{N}\ast f(z)\right\vert ^{r}dz \\
&\lesssim &\text{ }|k|^{-m}N^{n-m}\left( l(P)\right) ^{-m}\int_{|z-x|\leq 2%
\sqrt{n}|k|l(P)}\left\vert \omega _{N}\ast f(z)\right\vert ^{r}\chi
_{P+kl(P)}(z)dz \\
&\lesssim &\text{ }|k|^{n-m}\left( Nl(P)\right) ^{n-m}\mathcal{M}\Big(%
|\omega _{N}\ast f|^{r}\chi _{P+kl(P)}\Big)(x).
\end{eqnarray*}%
Hence the left-hand side of $\mathrm{\eqref{exp3}}$ is bounded by 
\begin{eqnarray*}
&&\left( Nl(P)\right) ^{n-m}\left\Vert C\text{ }\mathcal{M}\Big(%
|k|^{-c_{\log }(\frac{1}{\tau })r}\frac{|\omega _{N}\ast f|^{r}\chi
_{P+kl(P)}}{\left\Vert \chi _{P}\right\Vert _{\tau (\cdot )}^{r}}\Big)%
\right\Vert _{p(\cdot )/r} \\
&\lesssim &\left( Nl(P)\right) ^{n-m}|k|^{-c_{\log }(\frac{1}{\tau }%
)r}\left\Vert \frac{|\omega _{N}\ast f|^{r}\chi _{P+kl(P)}}{\left\Vert \chi
_{P}\right\Vert _{\tau (\cdot )}^{r}}\right\Vert _{p(\cdot )/r},
\end{eqnarray*}%
after using the fact that $\mathcal{M}:L^{p(\cdot )/r}\rightarrow L^{p(\cdot
)/r}$ is bounded. By Lemma \ref{n-cube-est} (i), $\frac{\left\Vert \chi
_{P+kl(P)}\right\Vert _{\tau (\cdot )}}{\left\Vert \chi _{P}\right\Vert
_{\tau (\cdot )}}\leq c|k|^{c_{\log }(\frac{1}{\tau })},$ with $c>0$
independent of $N,h$ and $k$. Hence the last expression is bounded by%
\begin{equation*}
c\left( Nl(P)\right) ^{n-m}\left\Vert \frac{|\omega _{N}\ast f|^{r}\chi
_{P+kl(P)}}{\left\Vert \chi _{P+kl(P)}\right\Vert _{\tau (\cdot )}^{r}}%
\right\Vert _{p(\cdot )/r}.
\end{equation*}%
Since $m$ can be taken large enough such that $m>2n+c_{\log }(\frac{1}{\tau }%
)r$, then the second sum in $\mathrm{\eqref{sec-sum}}$ is bounded by%
\begin{eqnarray*}
&&\left( Nl(P)\right) ^{n-m}\sum_{k\in \mathbb{Z}^{n},|k|>4\sqrt{n}%
}|k|^{n+c_{\log }(\frac{1}{\tau })r-m}\left\Vert \frac{|\omega _{N}\ast
f|^{r}}{\left\Vert \chi _{P+kl(P)}\right\Vert _{\tau (\cdot )}^{r}}\chi
_{P+kl(P)}\right\Vert _{p(\cdot )/r} \\
&\lesssim &\left( Nl(P)\right) ^{n-m}\sum_{k\in \mathbb{Z}^{n},|k|>4\sqrt{n}%
}|k|^{n+c_{\log }(\frac{1}{\tau })r-m}\left\Vert \omega _{N}\ast
f\right\Vert _{L_{\tau (\cdot )}^{p(\cdot )}}^{r} \\
&\lesssim &\left( Nl(P)\right) ^{n-m}\left\Vert \omega _{N}\ast f\right\Vert
_{L_{\tau (\cdot )}^{p(\cdot )}}^{r}.
\end{eqnarray*}%
Hence the proof is complete. }

{\emph{Proof of Lemma }\ref{Key-lemma}. By similarity, we only consider ${%
\ell ^{\tau (\cdot ),q(\cdot )}(L^{p(\cdot )})}$ spaces. Let $P\in \mathcal{Q%
}$. In view of the proof of Lemma \ref{Alm-Hastolemma1} the problem can be
reduced to the case when $\ell ^{q(\cdot )}(L^{p(\cdot )})$ is a normed
space. Then%
\begin{eqnarray}
&&\left\Vert \left( \frac{g_{v}}{\left\Vert \chi _{P}\right\Vert _{\tau
(\cdot )}}\chi _{P}\right) _{v\geq v_{P}^{+}}\right\Vert _{\ell ^{q(\cdot
)}(L^{p(\cdot )})}  \label{est3} \\
&\leq &\left\Vert \left( \sum_{k=0}^{v_{P}^{+}}\frac{2^{-|k-v|\delta }f_{k}}{%
\left\Vert \chi _{P}\right\Vert _{\tau (\cdot )}}\chi _{P}\right) _{v\geq
v_{P}^{+}}\right\Vert _{\ell ^{q(\cdot )}(L^{p(\cdot )})}+\left\Vert \left(
\sum_{k=v_{P}^{+}}^{v}\cdot \cdot \cdot \right) _{v\geq
v_{P}^{+}}\right\Vert _{\ell ^{q(\cdot )}(L^{p(\cdot )})}  \notag \\
&&+\left\Vert \left( \sum_{k=v}^{\infty }\cdot \cdot \cdot \right) _{v\geq
v_{P}^{+}}\right\Vert _{\ell ^{q(\cdot )}(L^{p(\cdot )})}.  \notag
\end{eqnarray}%
The first norm is bounded by%
\begin{equation*}
\sum_{k=0}^{v_{P}^{+}}2^{(k-v_{P}^{+})\delta }\left\Vert \left( \frac{%
2^{(v_{P}^{+}-v)\delta }f_{k}}{\left\Vert \chi _{P}\right\Vert _{\tau (\cdot
)}}\chi _{P}\right) _{v\geq v_{P}^{+}}\right\Vert _{\ell ^{q(\cdot
)}(L^{p(\cdot )})}.
\end{equation*}%
Let $Q_{k,h}$ be a dyadic cube such that $P\subset Q_{k,h}$. Obviously $%
v_{Q_{k,h}}^{+}=k$ and by Lemma \ref{n-cube-est} we have $\frac{\left\Vert
\chi _{Q_{k,h}}\right\Vert _{\tau (\cdot )}}{\left\Vert \chi _{P}\right\Vert
_{\tau (\cdot )}}\lesssim 2^{n(v_{P}^{+}-k)/\tau ^{-}}$. Therefore the last
sum is bounded by%
\begin{eqnarray*}
&&\sum_{k=0}^{v_{P}^{+}}2^{(k-v_{P}^{+})(\delta -n/\tau ^{-})}\left\Vert
\left( \frac{f_{j}}{\left\Vert \chi _{Q_{k,h}}\right\Vert _{\tau (\cdot )}}%
\chi _{Q_{k,h}}\right) _{j\geq v_{Q_{k,h}}^{+}}\right\Vert _{\ell ^{q(\cdot
)}(L^{p(\cdot )})} \\
&\leq &\sum_{k=0}^{v_{P}^{+}}2^{(k-v_{P}^{+})(\delta -n/\tau
^{-})}\left\Vert (f_{v})_{v}\right\Vert _{\ell ^{\tau (\cdot ),q(\cdot
)}(L^{p(\cdot )})}\lesssim \left\Vert (f_{v})_{v}\right\Vert _{\ell ^{\tau
(\cdot ),q(\cdot )}(L^{p(\cdot )})},
\end{eqnarray*}%
since $\delta >n/\tau ^{-}$. Let $\sigma >\max (q^{+},\frac{q^{+}}{p^{-}})$
and $\left\Vert (f_{v})_{v}\right\Vert _{\ell ^{\tau (\cdot ),q(\cdot
)}(L^{p(\cdot )})}=1$. Then%
\begin{eqnarray*}
&&\sum\limits_{v=v_{P}^{+}}^{\infty }\left\Vert \left\vert \frac{%
\sum_{k=v_{P}^{+}}^{v}2^{(k-v)\delta }f_{v}}{\left\Vert \chi _{P}\right\Vert
_{\tau (\cdot )}}\right\vert ^{q\left( \cdot \right) }\chi _{P}\right\Vert _{%
\frac{p(\cdot )}{q\left( \cdot \right) }} \\
&=&\sum\limits_{v=v_{P}^{+}}^{\infty }\left\Vert \left\vert \frac{%
\sum_{k=v_{P}^{+}}^{v}2^{(k-v)\delta }f_{v}}{\left\Vert \chi _{P}\right\Vert
_{\tau (\cdot )}}\right\vert ^{q\left( \cdot \right) /\sigma }\chi
_{P}\right\Vert _{\frac{\sigma p(\cdot )}{q\left( \cdot \right) }}^{\sigma }
\\
&\leq &\sum\limits_{v=v_{P}^{+}}^{\infty }\left( \sum_{k=v_{P}^{+}}^{v}2^{%
\frac{(k-v)\delta q^{-}}{\sigma }}\left\Vert \left\vert \frac{f_{v}}{%
\left\Vert \chi _{P}\right\Vert _{\tau (\cdot )}}\right\vert ^{q\left( \cdot
\right) /\sigma }\chi _{P}\right\Vert _{\frac{\sigma p(\cdot )}{q\left(
\cdot \right) }}\right) ^{\sigma } \\
&\lesssim &\sum\limits_{v=v_{P}^{+}}^{\infty }\left\Vert \left\vert \frac{%
f_{v}}{\left\Vert \chi _{P}\right\Vert _{\tau (\cdot )}}\right\vert
^{q\left( \cdot \right) }\chi _{P}\right\Vert _{\frac{p(\cdot )}{q\left(
\cdot \right) }}\leq 1,
\end{eqnarray*}%
by Lemma \ref{lq-inequality}. The desired estimate is completed by the
scaling argument. Now the last norm in $\mathrm{\eqref{est3}}$ is bounded by%
\begin{eqnarray*}
&&\left\Vert \left( \sum_{i=0}^{\infty }\frac{2^{-i\delta }f_{i+v}}{%
\left\Vert \chi _{P}\right\Vert _{\tau (\cdot )}}\chi _{P}\right) _{v\geq
v_{P}^{+}}\right\Vert _{\ell ^{q(\cdot )}(L^{p(\cdot )})} \\
&\leq &\sum_{i=0}^{\infty }2^{-i\delta }\left\Vert \left( \frac{f_{k}}{%
\left\Vert \chi _{P}\right\Vert _{\tau (\cdot )}}\chi _{P}\right) _{k\geq
v_{P}^{+}+i}\right\Vert _{\ell ^{q(\cdot )}(L^{p(\cdot )})} \\
&\leq &\sum_{i=0}^{\infty }2^{-i\delta }\left\Vert \left( \frac{f_{k}}{%
\left\Vert \chi _{P}\right\Vert _{\tau (\cdot )}}\chi _{P}\right) _{k\geq
v_{P}^{+}}\right\Vert _{\ell ^{q(\cdot )}(L^{p(\cdot )})}\lesssim \left\Vert
(f_{v})_{v}\right\Vert _{\ell ^{\tau (\cdot ),q(\cdot )}(L^{p(\cdot )})}.
\end{eqnarray*}%
Hence the lemma is proved. }

{\emph{Proof of Lemma }\ref{lamda-equi}. First we consider the space $%
b_{p(\cdot ),q(\cdot )}^{\alpha (\cdot ),\tau (\cdot )}$. Obviously, $%
\left\Vert \lambda \right\Vert _{b_{p\left( \cdot \right) ,q\left( \cdot
\right) }^{\alpha \left( \cdot \right) ,\tau (\cdot )}}\leq \left\Vert
\lambda _{r,d}^{\ast }\right\Vert _{b_{p\left( \cdot \right) ,q\left( \cdot
\right) }^{\alpha \left( \cdot \right) ,\tau (\cdot )}}$. Let us prove that $%
\left\Vert \lambda _{r,d}^{\ast }\right\Vert _{b_{p\left( \cdot \right)
,q\left( \cdot \right) }^{\alpha \left( \cdot \right) ,\tau (\cdot )}}\leq
c\left\Vert \lambda \right\Vert _{b_{p\left( \cdot \right) ,q\left( \cdot
\right) }^{\alpha \left( \cdot \right) ,\tau (\cdot )}}$. By the scaling
argument, it suffices to consider the case $\left\Vert \lambda \right\Vert
_{b_{p\left( \cdot \right) ,q\left( \cdot \right) }^{\alpha \left( \cdot
\right) ,\tau (\cdot )}}=1$ and show that the modular of a constant times
the sequence on the left-hand side is bounded. It suffices to prove that%
\begin{eqnarray}
&&\left\Vert \left\vert \frac{c\sum\limits_{m\in \mathbb{Z}^{n}}2^{v(\alpha
\left( \cdot \right) +n/2)}\lambda _{v,m,r,d}^{\ast }\chi _{v,m}}{\left\Vert
\chi _{P}\right\Vert _{\tau (\cdot )}}\right\vert ^{q(\cdot )}\chi
_{P}\right\Vert _{\frac{p(\cdot )}{q(\cdot )}}  \notag \\
&\leq &\left\Vert \left\vert \sum_{i=0}^{\infty }2^{\epsilon i/r}\frac{%
\sum\limits_{m\in \mathbb{Z}^{n}}2^{v(\alpha \left( \cdot \right)
+n/2)}\lambda _{v,m}\chi _{v,m}}{\left\Vert \chi
_{Q(c_{P},2^{i-v_{P}})}\right\Vert _{\tau (\cdot )}}\right\vert ^{q(\cdot
)}\chi _{Q(c_{P},2^{i-v_{P}})}\right\Vert _{\frac{p(\cdot )}{q(\cdot )}%
}+2^{-v}=\delta  \label{delta-est}
\end{eqnarray}%
where, $v\geq v_{P}^{+},$ $\epsilon =(n-d+a+n/\tau ^{-})/2$, $P\in \mathcal{Q%
}$ and $Q(c_{P},2^{i-v_{P}})$ is the cube concentric with $P$ having the
side length $2^{i-v_{P}}$. Therefore,%
\begin{equation*}
\sum_{v=v_{P}^{+}}^{\infty }\left\Vert \left\vert \tfrac{\sum\limits_{m\in 
\mathbb{Z}^{n}}2^{v(\alpha \left( \cdot \right) +n/2)}\lambda
_{v,m,r,d}^{\ast }\chi _{v,m}}{\left\Vert \chi _{P}\right\Vert _{\tau (\cdot
)}}\right\vert ^{q(\cdot )}\chi _{P}\right\Vert _{\frac{p(\cdot )}{q(\cdot )}%
}\lesssim 1
\end{equation*}%
for any dyadic cube $P$. The claim can be reformulated as showing that%
\begin{equation*}
\left\Vert \delta ^{-1}\left\vert \frac{c\sum\limits_{m\in \mathbb{Z}%
^{n}}2^{v(\alpha \left( \cdot \right) +n/2)}\lambda _{v,m,r,d}^{\ast }\chi
_{v,m}}{\left\Vert \chi _{P}\right\Vert _{\tau (\cdot )}}\right\vert
^{q(\cdot )}\chi _{P}\right\Vert _{\frac{p(\cdot )}{q\left( \cdot \right) }%
}\leq 1,
\end{equation*}%
which is equivalent to%
\begin{equation}
\left\Vert \delta ^{-\frac{1}{q\left( \cdot \right) }}\frac{%
\sum\limits_{m\in \mathbb{Z}^{n}}2^{v(\alpha \left( \cdot \right)
+n/2)}\lambda _{v,m,r,d}^{\ast }\chi _{v,m}}{\left\Vert \chi _{P}\right\Vert
_{\tau (\cdot )}}\chi _{P}\right\Vert _{p(\cdot )}\leq c.  \label{key-est11}
\end{equation}%
For each $k\in \mathbb{N}_{0}$ we define $\Omega _{k}:=\{h\in \mathbb{Z}%
^{n}:2^{k-1}<2^{v}\left\vert 2^{-v}h-2^{-v}m\right\vert \leq 2^{k}\}$\ and $%
\Omega _{0}:=\{h\in \mathbb{Z}^{n}:2^{v}\left\vert
2^{-v}h-2^{-v}m\right\vert \leq 1\}$. Then for any $x\in Q_{v,m}\cap P$, $%
\sum_{h\in \mathbb{Z}^{n}}\frac{\delta ^{-\frac{r}{q\left( x\right) }%
}2^{vr\alpha \left( x\right) }|\lambda _{v,h}|^{r}}{%
(1+2^{v}|2^{-v}h-2^{-v}m|)^{d}}$ can be rewritten as%
\begin{eqnarray}
&&\sum\limits_{k=0}^{\infty }\sum\limits_{h\in \Omega _{k}}\frac{\delta ^{-%
\frac{r}{q\left( x\right) }}2^{vr\alpha \left( x\right) }\left\vert \lambda
_{v,h}\right\vert ^{r}}{\left( 1+2^{v}\left\vert 2^{-v}h-2^{-v}m\right\vert
\right) ^{d}}  \notag \\
&\lesssim &\sum\limits_{k=0}^{\infty }2^{-dk}\sum\limits_{h\in \Omega
_{k}}\delta ^{-\frac{r}{q\left( x\right) }}2^{vr\alpha \left( x\right)
}\left\vert \lambda _{v,h}\right\vert ^{r}  \notag \\
&=&\sum\limits_{k=0}^{\infty }2^{(n-d)k+(v-k)n+vr\alpha \left( x\right)
}\delta ^{-\frac{r}{q\left( x\right) }}\int\limits_{\cup _{z\in \Omega
_{k}}Q_{v,z}}\sum\limits_{h\in \Omega _{k}}\left\vert \lambda
_{v,h}\right\vert ^{r}\chi _{v,h}(y)dy.  \label{key-est2}
\end{eqnarray}%
Let $x\in Q_{v,m}\cap P$ and $y\in \cup _{z\in \Omega _{k}}Q_{v,z}$, then $%
y\in Q_{v,z}$ for some $z\in \Omega _{k}$ and $2^{k-1}<2^{v}\left\vert
2^{-v}z-2^{-v}m\right\vert \leq 2^{k}$. From this it follows that $y$ is
located in some cube $Q(x,2^{k-v+3})$. In addition, from the fact that%
\begin{equation*}
\left\vert y_{i}-(c_{P})_{i}\right\vert \leq \left\vert
y_{i}-x_{i}\right\vert +\left\vert x_{i}-(c_{P})_{i}\right\vert \leq
2^{k-v+2}+2^{-v_{P}-1}<2^{k-v_{P}+3},\text{ \ }i=1,...,n,
\end{equation*}%
we have $y$ is located in some cube $Q(c_{P},2^{k-v_{P}+4})$. Since $1/q$ is
log-H\"{o}lder continuous and $\delta \in \lbrack 2^{-v},1+2^{-v}]$, we have 
}

\begin{equation*}
\delta ^{\frac{1}{q(x)}-\frac{1}{q(y)}}\leq 2^{\left\vert \frac{1}{q(x)}-%
\frac{1}{q(y)}\right\vert (2v+1)}\leq 2^{\frac{c_{\log }(q)\text{ }(2v+1)}{%
\log (e+\frac{1}{\left\vert x-y\right\vert })}}\leq 2^{\frac{c_{\log }(q)%
\text{ }(2v+1)}{v-k-h_{n}}}\lesssim \text{ }2^{2c_{\log }(q)k}
\end{equation*}%
for any $k<\max (0,v-h_{n})$ and any $y\in Q(x,2^{k-v+3})$, with $h_{n}\in 
\mathbb{N}$. If $k\geq \max (0,v-h_{n})$ then since again $\delta \in
\lbrack 2^{-v},1+2^{-v}]$, $\delta ^{\frac{1}{q(x)}-\frac{1}{q(y)}}\leq c$ $%
2^{\left\vert \frac{1}{q(x)}-\frac{1}{q(y)}\right\vert (2v+1)}\leq c$ $2^{2(%
\frac{1}{q^{-}}-\frac{1}{q^{+}})k}$. Also since $\alpha $ is log-H\"{o}lder
continuous we can prove that 
\begin{equation*}
2^{v(\alpha \left( x\right) -\alpha \left( y\right) )}\lesssim \left\{ 
\begin{array}{ccc}
2^{c_{\log }(\alpha )k} & \text{if} & k<\max (0,v-h_{n}) \\ 
2^{(\alpha ^{+}-\alpha ^{-})k} & \text{if} & k\geq \max (0,v-h_{n}),%
\end{array}%
\right.
\end{equation*}%
where $c>0$ not depending on $v$ and $k$. Therefore, $\mathrm{%
\eqref{key-est2}}$ does not exceed%
\begin{eqnarray*}
&&c\sum\limits_{k=0}^{\infty
}2^{(n-d+a)k+(v-k)n}\int\limits_{Q(x,2^{k-v+3})}\delta ^{-\frac{r}{q\left(
y\right) }}2^{v\alpha \left( y\right) r}\sum\limits_{h\in \Omega
_{k}}\left\vert \lambda _{v,h}\right\vert ^{r}\chi _{v,h}(y)\chi
_{Q(c_{P},2^{k-v_{P}+4})}dy \\
&\lesssim &\sum\limits_{k=0}^{\infty }2^{(n-d+a)k}\mathcal{M}\Big(%
\sum\limits_{h\in \Omega _{k}}\delta ^{-\frac{r}{q\left( \cdot \right) }%
}2^{v\alpha \left( \cdot \right) r}\left\vert \lambda _{v,h}\right\vert
^{r}\chi _{v,h}\chi _{Q(c_{P},2^{k-v_{P}+4})}\Big)(x).
\end{eqnarray*}%
Hence the left-hand side of $\mathrm{\eqref{key-est11}}$ is bounded by 
\begin{eqnarray*}
&&c\left\Vert \frac{\sum\limits_{k=0}^{\infty }2^{(n-d+a)k}\mathcal{M}\Big(%
\sum\limits_{h\in \Omega _{k}}\delta ^{-\frac{r}{q\left( \cdot \right) }%
}2^{v(\alpha \left( \cdot \right) +n/2)r}\left\vert \lambda
_{v,h}\right\vert ^{r}\chi _{v,h}\chi _{Q(c_{P},2^{k-v_{P}+4})}\Big)}{%
\left\Vert \chi _{P}\right\Vert _{\tau (\cdot )}^{r}}\right\Vert _{p(\cdot
)/r}^{1/r} \\
&\lesssim &\left( \sum\limits_{k=0}^{\infty }2^{\epsilon k}\left\Vert
\sum\limits_{i=0}^{\infty }2^{\epsilon i/r}\frac{\sum\limits_{h\in \Omega
_{i}}\delta ^{-\frac{1}{q\left( \cdot \right) }}2^{v(\alpha \left( \cdot
\right) +n/2)}\left\vert \lambda _{v,h}\right\vert \chi _{v,h}\chi
_{Q(c_{P},2^{i-v_{P}})}}{\left\Vert \chi _{Q(c_{P},2^{i-v_{P}})}\right\Vert
_{\tau (\cdot )}}\right\Vert _{p(\cdot )}^{r}\right) ^{1/r} \\
&\lesssim &\left( \sum\limits_{k=0}^{\infty }2^{\epsilon k}\right) ^{1/r},
\end{eqnarray*}%
where on the first estimate we use Lemma \ref{n-cube-est} and the
boundedness of the maximal function on $L^{p/r}$ (since $r<p^{-}$), and for
the last estimate we use the fact that 
\begin{equation*}
\left\Vert \sum\limits_{i=0}^{\infty }2^{\epsilon i/r}\frac{%
\sum\limits_{h\in \mathbb{Z}^{n}}\delta ^{-\frac{1}{q\left( \cdot \right) }%
}2^{v(\alpha \left( \cdot \right) +n/2)}\left\vert \lambda _{v,h}\right\vert
\chi _{v,h}}{\left\Vert \chi _{Q(c_{P},2^{i-v_{P}})}\right\Vert _{\tau
(\cdot )}}\chi _{Q(c_{P},2^{i-v_{P}})}\right\Vert _{p(\cdot )}\lesssim 1
\end{equation*}%
and $d$ sufficiently large such that $d>n+a+n/\tau ^{-}$.

Now we consider the space\ $\widetilde{b}_{p(\cdot ),q(\cdot )}^{\alpha
(\cdot ),p(\cdot )}$. Obviously we need only to prove $\mathrm{%
\eqref{delta-est}}$\textrm{\ }with $\left\vert
Q(c_{P},2^{i-v_{P}})\right\vert ^{1/p(\cdot )}$ in place of $\left\Vert \chi
_{Q(c_{P},2^{i-v_{P}})}\right\Vert _{\tau (\cdot )}$, $\epsilon
=(n-d+a+c_{\log }(1/p)+n/p^{-})/2$, $P\in \mathcal{Q}$, $\left\vert
P\right\vert \leq 1$ ($\left\Vert \chi _{P}\right\Vert _{p(\cdot )}$ in
place of $\left\Vert \chi _{P}\right\Vert _{\tau (\cdot )}$). We use the
same arguments above, we obtain that the left-hand side of $\mathrm{%
\eqref{key-est11}}$ (with power $r$) is bounded by%
\begin{equation}
c\sum\limits_{k=0}^{\infty }2^{(n-d+a)k}\left\Vert \frac{1}{\left\Vert \chi
_{P}\right\Vert _{p(\cdot )}^{r}}\sum\limits_{h\in \Omega _{k}}\delta ^{-%
\frac{r}{q\left( \cdot \right) }}2^{v(\alpha \left( \cdot \right)
+n/2)r}\left\vert \lambda _{v,h}\right\vert ^{r}\chi _{v,h}\chi
_{Q(c_{P},2^{k-v_{P}+4})}\right\Vert _{p(\cdot )/r}  \label{delta-est1}
\end{equation}%
for any dyadic cube $P\in \mathcal{Q}$, with $|P|\leq 1$. Observe that 
\begin{eqnarray*}
2^{nv_{P}/p(x)} &\lesssim &\left( 1+2^{v_{P}}\left\vert x-y\right\vert
\right) ^{c_{\log }(1/p)}2^{nv_{P}/p(y)}\lesssim 2^{kc_{\log
}(1/p)}2^{nv_{P}/p(y)} \\
&\lesssim &2^{k\left( c_{\log }(1/p)+n/p^{-}\right) }\left\vert
Q(c_{P},2^{k-v_{P}+4})\right\vert ^{-1/p(y)}
\end{eqnarray*}%
for any $x\in P$ and any $y\in Q(c_{P},2^{k-v_{P}+4})$. Hence, $\mathrm{%
\eqref{delta-est1}}$ is bounded by%
\begin{equation*}
c\sum\limits_{k=0}^{\infty }2^{2\epsilon k}\left\Vert \frac{%
\sum\limits_{h\in \Omega _{k}}\delta ^{-\frac{r}{q\left( \cdot \right) }%
}2^{v(\alpha \left( \cdot \right) +n/2)r}\left\vert \lambda
_{v,h}\right\vert ^{r}\chi _{v,h}\chi _{Q(c_{P},2^{k-v_{P}+4})}}{\left\vert
Q(c_{P},2^{k-v_{P}+4})\right\vert ^{r/p(\cdot )}}\right\Vert _{p(\cdot
)/r}\lesssim 1.
\end{equation*}%
The proof of the lemma is thus complete.

{\emph{The convergence of }\eqref{new-rep}. Let $\varphi \in \mathcal{S(}%
\mathbb{R}^{n})$. By $\mathrm{\eqref{supp-cond}}$\textrm{-}$\mathrm{%
\eqref{diff-cond}}$\textrm{-}$\mathrm{\eqref{mom-cond}}$ and the Taylor
expansion of $\varphi $ up to order $L$\ with respect to the off-points $%
x_{Q_{v,m}}$, we obtain for fixed $v$%
\begin{eqnarray*}
&&\int_{\mathbb{R}^{n}}\sum\limits_{m\in \mathbb{Z}^{n}}\lambda
_{v,m}\varrho _{v,m}(y)\varphi (y)dy \\
&=&\int_{\mathbb{R}^{n}}\sum\limits_{m\in \mathbb{Z}^{n}}\lambda
_{v,m}\varrho _{v,m}(y)\left( \varphi (y)-\sum\limits_{\left\vert \beta
\right\vert \leq L}(y-x_{Q_{v,m}})^{\beta }\frac{\partial ^{\alpha }\varphi
(x_{Q_{v,m}})}{\beta !}\right) dy.
\end{eqnarray*}%
The last factor in the integral can be uniformly estimated from the above by%
\begin{equation*}
c\text{ }2^{-v(L+1)}(1+\left\vert y\right\vert ^{2})^{-M/2}\sup_{x\in 
\mathbb{R}^{n}}(1+\left\vert x\right\vert ^{2})^{M/2}\sum\limits_{\left\vert
\beta \right\vert \leq L+1}\left\vert \partial ^{\alpha }\varphi
(x)\right\vert ,
\end{equation*}%
where $M>0$ is at our disposal. Let $0<t<\left( p(\cdot )\left( 1-\frac{1}{%
\min (1,p^{-})}\right) \right) ^{-}+1$ and $s(x)=\alpha (x)+\frac{n}{p(x)}%
(t-1)$ be such that $L+1>-\alpha (\cdot )+n\left( \frac{1}{\min (1,p^{-})}%
-1\right) >-s(\cdot )$. Since $\varrho _{v,m}$ are $\left[ K,L\right] $%
-atoms, then for every $S>0$, we have $\left\vert \varrho
_{v,m}(y)\right\vert \leq c2^{vn/2}\left( 1+2^{v}\left\vert
y-x_{Q_{v,m}}\right\vert \right) ^{-S}$. Therefore,%
\begin{eqnarray*}
&&\left\vert \int_{\mathbb{R}^{n}}\sum\limits_{m\in \mathbb{Z}^{n}}\lambda
_{v,m}\varrho _{v,m}(y)\varphi (y)dy\right\vert  \\
&\leq &c\text{ }2^{-v(L+1)}\int_{\mathbb{R}^{n}}\sum\limits_{m\in \mathbb{Z}%
^{n}}2^{vn/2}\left\vert \lambda _{v,m}\right\vert \frac{(1+\left\vert
y\right\vert ^{2})^{-M/2}}{\left( 1+2^{v}\left\vert y-x_{Q_{v,m}}\right\vert
\right) ^{S}}dy \\
&=&2^{-v(L+1)}\sum\limits_{h\in \mathbb{Z}^{n}}\int_{Q_{0,h}}\cdot \cdot
\cdot dy.
\end{eqnarray*}%
Applying Lemma \ref{Conv-est1} to obtain%
\begin{equation*}
\sum\limits_{m\in \mathbb{Z}^{n}}\left\vert \lambda _{v,m}\right\vert \left(
1+2^{v}\left\vert y-x_{Q_{v,m}}\right\vert \right) ^{-S}\lesssim
\sum\limits_{m\in \mathbb{Z}^{n}}\left\vert \lambda _{v,m}\right\vert \eta
_{v,S}\ast \chi _{v,m}(y).
\end{equation*}%
We split $M$ into $R+S$. Since we have in addition the factor $(1+\left\vert
y\right\vert ^{2})^{-S/2}$, H\"{o}lder's inequality, the fact that $%
\left\Vert \chi _{Q_{0,h}}\right\Vert _{\tau (\cdot )}\approx \left\Vert
\chi _{Q_{0,h}}\right\Vert _{(p(\cdot )/t)^{\prime }}\approx 1$, see $%
\mathrm{\eqref{DHHR2}}$\textrm{, }and $(1+\left\vert y\right\vert
^{2})^{-R/2}\lesssim (1+\left\vert h\right\vert ^{2})^{-R/2}$ give that the
term $\left\vert \int_{\mathbb{R}^{n}}\cdot \cdot \cdot dy\right\vert $ is
bounded by%
\begin{eqnarray*}
&&c\ 2^{-v(L+1)}\sum\limits_{h\in \mathbb{Z}^{n}}(1+\left\vert h\right\vert
^{2})^{-R/2}\left\Vert \frac{\eta _{v,S}\ast \left[ \sum\limits_{m\in 
\mathbb{Z}^{n}}2^{vn/2}\left\vert \lambda _{v,m}\right\vert \chi _{v,m}%
\right] }{\left\Vert \chi _{Q_{0,h}}\right\Vert _{\tau (\cdot )}}\chi
_{Q_{0,h}}\right\Vert _{p(\cdot )/t} \\
&\lesssim &\sup_{P\in \mathcal{Q},j\geq j_{P}^{+}}\left\Vert \frac{%
2^{(s(\cdot )+n/2)j}\sum\limits_{m\in \mathbb{Z}^{n}}\left\vert \lambda
_{j,m}\right\vert \chi _{j,m}}{\left\Vert \chi _{P}\right\Vert _{\tau (\cdot
)}}\chi _{P}\right\Vert _{p(\cdot )/t}\lesssim \left\Vert \lambda
\right\Vert _{b_{p(\cdot )/t,\infty }^{s(\cdot ),\tau (\cdot )}},
\end{eqnarray*}%
where the first inequality follows by Lemma \ref{Alm-Hastolemma1}, $%
L+1+s(\cdot )>0$ and by taking $R$ large enough. The convergence of $\mathrm{%
\eqref{new-rep}}$ is now clear by the embeddings $\left\Vert \lambda
\right\Vert _{b_{p(\cdot ),q(\cdot )}^{\alpha ,\tau (\cdot
)}}\hookrightarrow \left\Vert \lambda \right\Vert _{b_{p(\cdot ),\infty
}^{\alpha ,\tau (\cdot )}}\hookrightarrow \left\Vert \lambda \right\Vert
_{b_{p(\cdot )/t,\infty }^{s(\cdot ),\tau (\cdot )}}$. The proof is
completed.\vskip5pt }

\textbf{Acknowledgements}\vskip5pt

We thank the anonymous referees for pointing the references [7, 15, 18, 19,
22, 25, 29, 35-39, 43-44] out to us and for the valuable comments and
suggestions, especially on the proof of Lemma \ref{r-trick}\textbf{\ }and
the motivation of introducing the fourth parameter $\tau $.


\begin{thebibliography}{99}
\bibitem{AH} A. Almeida and P. H\"{a}st\"{o}, \emph{Besov spaces with
variable smoothness and integrability}, J. Funct. Anal. \textbf{258} (2010),
1628--1655.

\bibitem{M07} M. Bownik, \emph{Anisotropic Triebel-Lizorkin Spaces with
Doubling Measures}, The Journal of Geometric Analysis. \textbf{17} (2007),
no. 3, 337--424.

\bibitem{CFMP} D. Cruz-Uribe, A. Fiorenza, J. M, Martell and C. P\'{e}rez, 
\emph{The boundedness of classical operators in variable $L^{p}$\ spaces},
Ann. Acad. Sci. Fenn. Math. \textbf{13} (2006), 239--264.

\bibitem{Di} L. Diening, \emph{Maximal function on generalized Lebesque
spaces $L^{p(\cdot )}$}, Math. Inequal. Appl. \textbf{7} (2004), no. 2,
245--253.

\bibitem{DHR} L. Diening, P. H\"{a}st\"{o} and S. Roudenko,\emph{\ Function
spaces of variable smoothness and integrability}, J. Funct. Anal. \textbf{256%
} (2009),no. 6, 1731--1768.

\bibitem{DHHR} L. Diening, P. Harjulehto, P. H\"{a}st\"{o} and M. R\r{u}\v{z}%
i\v{c}ka,\emph{\ Lebesgue and Sobolev spaces with variable\ exponents},
Lecture Notes in Mathematics, vol. 2017, Springer-Verlag, Berlin 2011.

\bibitem{DHHMS} L. Diening, P. Harjulehto, P. H\"{a}st\"{o}, Y. Mizuta and
T. Shimomura, \emph{Maximal functions in variable exponent spaces: limiting
cases of the exponent}, Ann. Acad. Sci. Fenn. Math. \textbf{34} (2009), no.
2, 503--522.

\bibitem{D1} D. Drihem, S\emph{ome embeddings and equivalent norms of the $%
\mathcal{L}_{p,q}^{\lambda ,s}$\ spaces}, Funct. Approx. Comment. Math. 
\textbf{41} (2009), no. 1, 15--40.

\bibitem{D3} D. Drihem, \emph{Atomic decomposition of Besov spaces with
variable smoothness and integrability}, J. Math. Anal. Appl. \textbf{389}
(2012), no. 1, 15--31.

\bibitem{D4} D. Drihem, \emph{Atomic decomposition of Besov-type and
Triebel-Lizorkin-type spaces}, Sci. China. Math. \textbf{56} (2013), no. 5,
1073--1086.

\bibitem{D5} D. Drihem, \emph{Some properties of variable Besov-type spaces}%
, Funct. Approx. Comment. Math. \textbf{52} (2015), no. 2, 193-221.

\bibitem{FJ86} M. Frazier and B. Jawerth, \emph{Decomposition of Besov spaces%
}, Indiana Univ. Math. J. \textbf{34} (1985), 777--799.

\bibitem{FJ90} M. Frazier and B. Jawerth, \emph{A discrete transform and
decomposition of distribution spaces}, J. Funct. Anal. \textbf{93}  (1990),
no. 1, 34--170.

\bibitem{FX11} J. Fu and J. Xu, \emph{Characterizations of Morrey type Besov
and Triebel-Lizorkin spaces with variable exponents}, J. Math. Anal. Appl. 
\textbf{381} (2011), 280--298.

\bibitem{HN07} L. Hedberg and Y. Netrusov, \emph{An axiomatic approach to
function spaces, spectral synthesis, and Luzin approximation}, Mem. Amer.
Math. Soc. \textbf{188} (2007), vi+97 pp.

\bibitem{KV121} H. Kempka and J. Vyb\'{\i}ral, \emph{A note on the spaces of
variable integrability and summability of Almeida and H\"{a}st\"{o}}, Proc.
Amer. Math. Soc. \textbf{141} (2013), no. 9, 3207--3212.

\bibitem{KV122} H. Kempka and J. Vyb\'{\i}ral, \emph{Spaces of variable
smoothness and integrability: Characterizations by local means and ball
means of differences}, J. Fourier Anal. Appl. \textbf{18}, (2012), no. 4,
852--891.

\bibitem{KY95} H. Kozono and M. Yamazaki, \emph{Semilinear heat equations
and the Navier-Stokes equation with distributions in new function spaces as
initial data}, Comm. PDE. \textbf{19} (1994), 959--1014.

\bibitem{LYYUS13} Y. Liang, Y. Sawano, T. Ullrich, D. Yang and W. Yuan, 
\emph{A new framework for generalized Besov-type and Triebel-Lizorkin-type
spaces}, Dissertationes Math. (Rozprawy Mat.) \textbf{489} (2013).

\bibitem{Ma03} A. L, Mazzucato, \emph{Besov-Morrey spaces: function space
theory and applications to non-linear PDE}, Trans. Amer. Math. Soc. \textbf{%
355} (2003), no. 4, 1297--1364.

\bibitem{M38} C. B., Jr. Morrey, \emph{On the solutions of quasi-linear
elliptic partial differential equations}, Trans.\ Amer. Math. Soc. \textbf{43%
} (1938), no. 1, 126--166.

\bibitem{N84} Y.V. Netrusov, \emph{Some imbedding theorems for spaces of
Besov-Morrey type}. (Russian) Numerical methods and questions in the
organization of calculations, 7. Zap. Nauchn. Sem. Leningrad. Otdel. Mat.
Inst. Steklov. (LOMI) \textbf{139} (1984), 139--147.

\bibitem{P} J. Peetre, \emph{On the theory of $\mathcal{L}_{p,\lambda }$\
spaces}, J. Funct. Anal. \textbf{4} (1969), 71--87.

\bibitem{Ru00} M.\ R\r{u}\v{z}i\v{c}ka, \emph{Electrorheological fluids:
modeling and mathematical theory}, Lecture Notes in Mathematics, 1748,
Springer-Verlag, Berlin, 2000.

\bibitem{Ry01} V.S. Rychkov, \emph{Littlewood-Paley theory and function
spaces with $A_{p}^{\text{loc}}$ weights}, Math. Nachr. \textbf{224} (2001),
145--180.

\bibitem{Sa08} Y. Sawano,\ \emph{Wavelet characterizations of Besov-Morrey
and Triebel-Lizorkin-Morrey spaces}, Funct. Approx. Comment. Math. \textbf{38%
} (2008), 93--107.

\bibitem{Sa09} Y. Sawano, \emph{A note on Besov-Morrey spaces and
Triebel-Lizorkin-Morrey spaces}, Acta Math.\ Sinica, English Series. \textbf{%
25} (2009), no. 8, 1223--1242.

\bibitem{Sa10} Y. Sawano, \emph{Brezis-Gallou\"{e}t-Wainger type inequality
for Besov-Morrey spaces}, Studia Math. \textbf{196} (2010), 91--101.

\bibitem{SYY10} Y. Sawano, D. Yang and W. Yuan, \emph{New applications of
Besov-type and Triebel-Lizorkin-type spaces}, J. Math. Anal. Appl. \textbf{%
363} (2010), 73--85.

\bibitem{ST89} J.-O. Str\"{o}mberg and A. Torchinsky, \emph{Weighted Hardy
spaces}, Lecture Notes in Math. 1381, Springer, Berlin, 1989.

\bibitem{SW} E. M. Stein and G. Weiss, I\emph{ntroduction to Fourier
analysis on Euclidean spaces}, Princeton Univ. Press, Princeton, NJ, 1971.

\bibitem{T1} H. Triebel, \emph{Theory of Function Spaces}, Birkh\"{a}user
Verlag, Basel, 1983.

\bibitem{T2} H. Triebel, \emph{Theory of Function Spaces II}, Birkh\"{a}user
Verlag, Basel, 1992.

\bibitem{T3} H. Triebel, \emph{Fractals and spectra}, Birkh\"{a}user, Basel
1997.

\bibitem{T4} H. Triebel, \emph{Local function spaces, heat and Navier-Stokes
equations}, EMS Tracts in Mathematics, 20. European Mathematical Society
(EMS), Zurich, 2013. x+232 pp.

\bibitem{YY1} D. Yang and W. Yuan, \emph{A new class of function spaces
connecting Triebel-Lizorkin spaces and $Q$\ spaces}, J. Funct. Anal. \textbf{%
255} (2008), 2760--2809.

\bibitem{YY2} D. Yang and W. Yuan, \emph{New Besov-type spaces and
Triebel-Lizorkin-type spaces including $Q$\ spaces}, Math. Z. \textbf{265}
(2010), 451--480.

\bibitem{YY13} D. Yang and W. Yuan, \emph{Relations among Besov-Type spaces,
Triebel-Lizorkin-Type spaces and generalized Carleson measure spaces}, Appl.
Anal, \textbf{92} (2013), no. 3, 549--561.

\bibitem{YSY13} W. Yuan, W. Sickel and D. Yang, \emph{On the coincidence of
certain approaches to smoothness spaces related to Morrey spaces}, Math.
Nachr. \textbf{286 } (2013), no. 14-15, 1571--1584.

\bibitem{YHSY} W. Yuan, D. Haroske, L. Skrzypczak, D. Yang, \emph{Embedding
properties of Besov-type spaces}, Applicable Analysis. \textbf{94 } (2015),
no. 2, 318--340.

\bibitem{YZW15} D. Yang, C. Zhuo and W. Yuan, \emph{Besov-Type Spaces with
Variable Smoothness and Integrability}, arXiv:1503.04512.

\bibitem{WYY} W. Yuan, W. Sickel and D. Yang, \emph{Morrey and Campanato
meet Besov, Lizorkin and Triebel}, Lecture Notes in Mathematics, vol. 2005,
Springer-Verlag, Berlin 2010.

\bibitem{Xu08} J. Xu, \emph{Variable Besov and Triebel-Lizorkin spaces},
Ann. Acad. Sci. Fenn. Math. \textbf{33} (2008), 511--522.

\bibitem{Xu09} J. Xu, An atomic decomposition of variable Besov and
Triebel-Lizorkin spaces, Armen. J. Math. \textbf{2} (2009), no. 1, 1--12.
\end{thebibliography}
\end{document}